\documentclass[a4paper]{amsart}
 \usepackage{amsmath,amssymb,amsthm,latexsym}
 \usepackage[cp1251]{inputenc}
  \theoremstyle{definition}
  \newtheorem{ddd}{Definition}[section]
  
  \theoremstyle{plain}
  \newtheorem{ttt}[ddd]{Theorem}
  
  \newtheorem{ccc}[ddd]{Corollary}
  \theoremstyle{remark}

 \newcommand{\supp}{\mathrm{supp}}

\begin{document}
    \title{Partially Commutative Groups And Lie Algebras}
    \thanks{The work is partially supported by RFBR (project 18--01--00100)}
    \author{E.\,N.\,Poroshenko, E.\,I.\,Timoshenko}
    \begin{abstract}
      This is a survey of results on partially commutative groups and
     partially commutative algebras.
    \end{abstract}
  \date{}
  \maketitle
  \section*{Introduction}

 Many algebraic structures are defined by graphs. Partially commutative algebraic structures
 are some of them. Let
 $\mathfrak M$ be a variety of algebraic structures of a functional signature
 $\Sigma$ containing a binary operation
 $\circ$. The case of commutative operation is trivial so suppose that this is not so.

 In this survey, by a graph we mean an undirected graph without loops. Graphs will be denoted
 by gothic letters.

 Let
 $\Delta=(X,E)$ be a graph (possibly infinite), with the set of vertices
 $X=\{x_1,x_2,\ldots\}$ and the set of edges
 $E=\{(x_i, x_j)\}$. For a variety
 $\mathfrak M$ define a \emph{partially commutative structure}
 $C(\mathfrak M, \Delta)$ on this variety as follows
 \begin{equation}\label{pcs}
   C(\mathfrak{M}, \Delta) = \langle X; x_i\circ x_j = x_j \circ x_i,\,\, \text{ if }
   (x_i, x_j) \in E,  \mathfrak M \rangle.
 \end{equation}
 Partially commutative structures appear in different areas of mathematics, for example, in
 computer science and robotics. By now, the most results in partially commutative
 structures have obtained for so called free partially commutative groups. A free partially
 commutative group is a partially commutative structure in the variety of all groups such
 that it is defined by an undirected graph without loops. Last years,  much attention is
 made to researches of partially commutative groups in soluble and nilpotent varieties.
 Partially commutative associative and Lie algebras are studied as well. In this survey,
 papers on partially commutative groups and Lie algebras are discussed. There are so many
 results in free partially commutative groups that a specific survey is needed for them.
 So, the results on free partially commutative groups are not included to this survey.
 Some information on these results can be found in
\cite{Cha07,DRT19}.

 There are two sections in the survey. In Sec.~%
\ref{sec1}, results for partially commutative groups of varieties are discussed. This
 section contains four subsections. In Subsec.~%
\ref{ssec11}, algebraic properties of
 partially commutative metabelian groups are observed. The results on the structure of the
 groups, their centralizers,  annihilators, bases, subgroups, inclusions into matrix groups, automorphism groups, and centralizer dimensions are described. At the end of the
 subsection, a decomposition of a group into a direct product is discussed. Such
 decompositions are considered not only for free partially commutative groups but also for
 partially commutative groups of varieties containing the variety of nilpotent groups
 of degree
 $\leq 2$.

 Let
 $C$ be a structure. The set
 $Th(C)$ of all first-order sentences of a signature
 $\Sigma$ which are true on
 $C$ is called the \emph{elementary theory} of
 $C$. Structures
 $C_1$ and
 $C_2$ are elementary equivalent if
 $Th(C_1) = Th(C_2)$.

 The universal theory or the
 $\forall$-theory of a structure
 $C$ is a subset of
 $Th(C)$ consisting of all
 $\forall$-formulas which are true on
 $C$.  Structures
 $C_1$ and
 $C_2$ are existentially equivalent if their existential theories coincide.

 Results on elementary and universal theories theories of partially commutative metabelian
 groups are considered in Subsec.~%
\ref{ssec12}. Information on varieties and
 prevarieties generated by partially commutative metabelian groups and on equations in one
 variable is also presented in Subsec.~%
\ref{ssec12}. The most attention is paid to universal
 theories, in particular, conditions of coincidence of two theories. In Subsec.~%
\ref{ssec13}, results on the structure and the universal theory of a partially commutative
  nilpotent group are considered. In Subsec.~%
\ref{ssec14}, theorems on centralizers and annihilators in partially commutative
 pro-$p$-groups are discussed.

 Sec.~%
\ref{sec2} mainly presents results on partially commutative Lie algebras in some varieties.
 This section consists of two subsections.

 In Subsec.~%
\ref{ssec21} algebraic results on partially commutative Lie algebras are discussed. In this
 subsection the results on isomorphisms, bases, annihilators and centralizers of partially
 commutative Lie algebras of some varieties are collected.

 In Subsec.~%
\ref{ssec22}, logical questions for partially commutative Lie algebras are discussed.
 Those are questions on universal and elementary theories of partially commutative Lie
 algebras.

  In Sec.~%
\ref{sec2},  there is a parallel presentation of results for partially commutative and
  partially commutative metabelian Lie algebras. Moreover, a description of a linear basis
  is also considered in the case of partially commutative nilpotent Lie algebras.

  Results on partially commutative groups defined by infinite graphs are also discussed in
  Sec.~%
\ref{sec2}, since there are analogous theorems for partially commutative Lie algebras in
  this section.

 Sec.~%
\ref{sec1} is written by the second author while
 Sec.
\ref{sec2} is written by the first one. 
  \section{Partially commutative groups} \label{sec1}

 This section is a survey of results for partially commutative groups on soluble varieties.

 Researches of groups defined by generators and defining relations form a large field
 of algebra. This field is called combinatorial group theory. It has a specific collection
 of problems and methods.  A lot of these methods have analogues in algebraic topology.
 Studies in combinatorial group theory are made heavily since the second half of the 
 20th century. One of the class of objects studied in combinatorial group theory consists of
 groups whose generators are the vertices of some graphs.

 Let
 $\Delta = (X,E)$  be a graph. For any variety
 $\mathfrak{M}$ and any graph
 $\Delta$ the partially commutative group
 $G(\mathfrak M, \Delta)$ in the variety
 $\mathfrak M$ has a representation
 $$G(\mathfrak M, \Delta) = \langle X \,| \, x_i x_j = x_jx_i,
   \text{ if } (x_i,x_j) \in E; \mathfrak{M}\rangle.$$
 Although generally problems considered for free partially commutative groups and for
 partially commutative groups in soluble varieties are same, research methods differ
 significantly. Methods effectively used for researches of algebraic properties in
 partially commutative groups of soluble varieties are those using modules over group rings,
 splitting extensions, Fox derivatives, etc.

 Let
 $G$ be a group,
 $g, h \in G$. Then we use the following notation.
 $[g, h] = g^{-1}h^{-1}gh$,
 $G' = [G, G]$. The subgroup
 $G'$ is called the commutant of
 $G$. The variety of metabelian groups
 $\mathfrak A^2$ is given by the identity
 $[[y_1, y_2],[y_3, y_4]] =1$. It means that this variety consists of groups
 $G$ having an abelian normal subgroup
 $A$ (possibly trivial) such that the quotient group
 $G/A$ is commutative. Denote a partially commutative group
 $G(\mathfrak A^2,\Delta)$ by
 $M_\Delta$ for short. Let
 $\mathfrak N_c$ be a variety of nilpotent groups of nilpotence degree at most
 $c$. This variety consists of all groups satisfying the identity
 $v_{c+1}=1$, where
 $v_2 = [y_1, y_2], v_{c+1} = [v_c, y_{c+1}]$. For a graph
 $\Delta$ denote by
 $M_{c,\Delta}$ the partially commutative group defined by
 $\Delta$ in the variety
 $\mathfrak A^2 \cap \mathfrak N_c$ .

 Some results on partially commutative metabelian groups can be found in
\cite{Tim13book}.

\subsection{Algebraic properties of partially commutative metabelian groups} \label{ssec11}

\subsubsection*{Torsion}
 If
 $\mathfrak N_2 \subseteq \mathfrak M$ then the quotient group
 $G(\mathfrak M,\Delta)/G'(\mathfrak M,\Delta)$ has no elements of finite order.

 Note that the periodic part of a group
 $G(\mathfrak M,\Delta)$ can be non-trivial. This is so, for example, if
 $\mathfrak M$ is a variety of centrally metabelian groups and
 $\Delta$ is a completely disconnected graph with at least four vertices
 (see \cite{Gup69}). The following theorem implies that partially commutative metabelian
 group has no elements of finite order.

 \begin{ttt}
\cite{Tim10} A group
  $M_\Delta$ can be approximated by nilpotent torsion-free groups.
 \end{ttt}

\subsubsection*{Center} Let
 $\Delta = (X; E)$ be a graph,
 $Y$ non-empty subset of the set
 $X$. We use the following notation.
 \begin{equation} \label{orthogonal}
   Y^\bot = \{x \in X\,|\,(x,y) \in E \text{ for all } y\in Y\}.
 \end{equation}

 Denote by
 $\langle  Y \rangle$ the group generated by
 $Y$.

 The following theorem describes the center of a partially commutative metabelian
 group, the quotient group by the center, and the relation of the center and the  commutant.

 \begin{ttt}
\cite{Tim10} Let
   $\Delta = (X;E)$ be a graph. Then the following statements hold.\\
   1) If
   $V^\bot $ is non-empty then
   $\mathcal{Z}(M_\Delta)=\langle X^\bot \rangle$, otherwise the center of
   $M_\Delta$ is trivial.\\
   2) If a subgraph
   $\Gamma$ of
   $\Delta$ is generated by the set
   $X\backslash X^\bot$ then
   $M_{\Delta}/ \mathcal Z(M_{\Delta}) \cong M_{\Gamma}$.\\
   3) The intersection of the center
   $\mathcal{Z}(M_{\Delta})$ and the commutant
   $M'_\Delta$ is trivial.
 \end{ttt}

\subsubsection* {Centralizers}
 Let
 $\Delta = (X;E)$ be a graph. It is interesting to consider partially commutative
 groups of varieties
 $\mathfrak M$ such that
 $u, v \in X$ commute in the corresponding group if and only if
 $u$ and
 $v$ are adjacent. Suppose that
 $\mathfrak{M}$ contains
 $\mathfrak{N}_2$. It turns out that for
 $u,v \in X$ the commutator
 $[u, v]$ is equal to one in
 $G(\mathfrak{M}, \Delta)$ if and only if
 $(u,v) \in E$.

 In Sec.~%
\ref{sec1}, we denote by
 $\overline{G}$ the quotient group
 $G/G'$ and by
 $\overline{g}$ the image of
 $g \in G$ in the group
 $\overline{G}$ via the natural homomorphism
 $G \to \overline{G}.$

 Let
 $G$  be a matabelian  (non-abelian) group. Its commutant
 $G'$ is a non-trivial abelian group and
 $G$ acts on
 $G'$ by conjugations:
 $c \mapsto g^{-1}cg$, for
 $g \in G$ and
 $c \in G'$. Since the elements in
 $G'$ act identically
 $G'$ is a right module on the integral group ring
 $\mathbb{Z}[\overline{G}]$. Denote the action of
 $\overline{g}$ on
 $c \in G'$ by
 $c^{\overline{g}}$.
 For elements
 $\alpha =\pm\overline{g}_1 \pm \ldots \pm \overline{g}_m \in \mathbb Z[\overline{G}]$
 and
 $c \in G'$ we put
 $$c^\alpha = {c^{\pm \overline{g}_1}} \cdot \ldots \cdot {c^{\pm\overline{g}_m}}.$$

 The centralizers of elements
 $x_i \in X$ and the centralizers in the commutant
 $\mathcal{C}(g) = C(g) \cap M_\Delta'$ of elements
 $g\in M_\Delta$ are described in the following theorem.

 \begin{ttt}
\cite{GT09,Tim10}  Let
   $X=\{x_1,\ldots,x_n\}$ be the set of vertices of the defining graph
   $\Delta$ of a group
   $M_\Delta$ and
   $\{x_1\}^\bot = \{x_2,\ldots, x_m\}$. The following statements hold.\\
   1) An element
   $g \in M_{\Delta}$ lies in the centralizer
   $C(x_1)$ of an element
   $x_1$ if and only if
   $$g= x_1^{l_1}\ldots x_m^{l_m}\prod_{2 \leq i < j \leq m} [x_i, x_j]^{\alpha_{ij}},$$
   where
   $l_1,\ldots, l_m \in \mathbb{Z}$,
   $\alpha_{ij} \in \mathbb{Z}[\overline{M_\Delta}]$. \\
   2) For any
   $m \leq n, \,\,1 \leq i_1 < \ldots < i_m \leq n$ and for any non-zero integers
   $q_1,\ldots, q_m$ the following equation holds
   $$\mathcal C (x_{i_1}^{q_1}\ldots x_{i_m}^{q_m}) =
     \mathcal C(x_{i_1}) \cap \ldots \cap \mathcal C(x_{i_m}). $$
 \end{ttt}

 Let us notice a couple of useful properties of centralizers of elements in groups
 $M_{\Delta}$ defined by trees or cycles. These properties are used to study the
 universal theories of partially commutative metabelian groups. In
\cite{GT09}, it was shown that the intersection of centralizers
 $\mathcal C(x_i)\cap \mathcal C (x_j)$ of two different elements
 $x_i, x_j \in X$ in
 $M_{\Delta}$ is trivial if
 $\Delta$ is a tree. If
 $\Delta$ is a cycle of length at least 4 then the intersection of centralizers
 $\mathcal C(x_i) \cap \mathcal{C} (x_j) \cap \mathcal{C}(x_l)$ of  tree different
 elements in
 $X$ is trivial
\cite{GT11-1}.

\subsubsection*{Annihilators}
  Let
  $c$ be an element in
  $M'_{\Delta}$. The \emph{annihilator}
  $Ann(c)$ of
  $c$ is an ideal of the ring
  $\mathbb{Z}[\overline{M}_{\Delta}]$, consisting of elements
  $\gamma$, such that
  $c^\gamma =1$.

  For any two non-adjacent vertices
  $x_i,x_j$ of
  $\Delta$ define the ideal
  $\mathcal{A}_{i,j}$ of the ring
  $\mathbb Z[\overline{M}_{\Delta}]$ as follows. If
  $x_i$ and
  $x_j$ lie in different connected components of
  $\Delta$ then put
  $\mathcal{A}_{i,j} =0$. Otherwise, let
  $a_i = x_i M'_{\Delta}$ and consider all paths with no returns
  $\{x_i, x_{i_1},\ldots, x_{i_m}, x_j\}$ connecting
  $x_i$ and
  $x_j$. To each path assign the element
  $(1-{a_{i_1}})\ldots (1-a_{x_{i_m}}) \in A = \mathbb
   Z[\overline{M}_{\Delta}] = \mathbb{Z}[a_1^{\pm 1},\dots,a_n^{\pm 1}]$. Let
  $\mathcal {A}_{i,j}$ be the ideal generated by these elements.

  \begin{ttt}
\cite{GT09} Let
   $\Delta = (X;E)$ be a graph with the set of vertices
   $X=\{x_1,\ldots, x_n\}$.  If
   $(x_i, x_j) \notin E$ then
   $Ann([x_i, x_j])= \mathcal{A}_{i,j}$.
 \end{ttt}

 The \emph{trivialization} of an element
 $\alpha$ in a group ring
 $\mathbb{Z}[G]$ is the image
 $\varepsilon(\alpha)$ of this element under the ring homomorphism
 $$\varepsilon: \mathbb{Z}[G] \to \mathbb{Z}$$
 extending the group homomorphism
 $G \to 1$. Useful properties of annihilators
 $\mathcal A_{ij}$ are given in the following theorem.

 \begin{ttt}\label{annihilators}
\cite{GT09} Let
   $X=\{x_1,\dots,x_n\}$ be the set of vertices of a graph
   $\Delta$ and
   $A=\mathbb{Z}[\overline{M}_{\Delta}] =
     \mathbb{Z}[a_1^{\pm1},\dots,a_n^{\pm 1}]$. Then the following statements hold.\\
   1) If
   $n \geq 3$,
   $a \in  A$, and
   $a(1-a_3)^2 \in \mathcal{A}_{1,2}$  then
   $a(1-a_3) \in \mathcal{A}_{1,2}$. \\
   2) If
   $n \geq 2$,
   $x_1,x_2$  are non-adjacent vertices of
   $\Delta$ and
   $a,\gamma \in A$ are such that
   $\varepsilon (\gamma) \neq 0$ and
   $a \gamma \in \mathcal{A}_{1,2}$ then
   $a \in \mathcal{A}_{1,2}$.\\
   3) if
   $n \geq 2$,
   $a \in   A$,
   $q,q_1,\dots,q_m$ are non-zero integers,
   $1<i_1<...<i_m \leq n$, and
   $$a(1-a_1^q a_{i_1}^{q_1}\dots a_{i_m}^{q_m})\in \mathcal{A}_{1,2}$$
   then
   $a \in\mathcal{A}_{1,2}.$\\
   4) if
   $n \geq 2$,
   $1\leq m \leq n$,
   $q_1,\dots,q_m$ are nonzero integers,
   $1 \leq i_1<\dots <i_m \leq n$, and
   $$a(1-a_{i_1}^{q_1}\dots a_{i_m}^{q_m})\in \mathcal{A}_{1,2}$$
   for an element
   $a \in A$ then all elements
   $a(1-a_{i_j})$ for
   $j=1,\dots, m$ are in
   $\mathcal{A}_{1,2}.$
 \end{ttt}

 Let us recall the definition of an associator. Consider a right module
 $L$ over a commutative ring
 $A$. A simple ideal
 $P$ of
 $A$ is associated with
 $L$ if there exists an element
 $0 \neq x \in L$ such that the annihilator of this element
 $$Ann(x)= \{a \in A\,|\,xa=0\}$$
 coincides with
 $P$. The set of ideals associated with the module
 $L$ is the \emph{associator} of
 $L$. In
\cite{GT11-1}, for a partially commutative group
 $M_\Delta$ the associator
 $\mathbb{Z}[\overline{M}_{\Delta}]$ of the module
 $M'_{\Delta}$ is described.

\subsubsection*{Basis and canonical representation of elements}
 The authors of
\cite{GT09} provided a theorem on canonical representation of elements of a partially
 commutative metabelian group. However, the proof of this theorem had a mistake and this
 was noticed in the paper
\cite{Tim10} of the second author of
 \cite{GT09}. Later on, in
\cite{ GT11-1,Tim10} theorems on a canonical representation of some elements in the
 commutatant of
 $M_\Delta$.  were proved. The presentation found there enabled to study the universal
 theory of a group
 $M_\Delta$ defined by a tree
 $\Delta$. However, a complete proof of a theorem on a canonical representation of
 elements in partially commutative metabelian group was given by the second author
 of this survey and
\cite{GT09} only in
 2020.

 The following theorem describes a basis of the commutant of a partially commutative
 metabelian group. It implies a canonical representation of elements of a the group.

 \begin{ttt}\cite{Tim20-1}\label{bas_comm}
   Let the set
   $X = \{x_1,\ldots,x_r\}$ of vertices of a graph
   $\Delta$ be ordered as follows
   $x_1 < x_2 < \ldots < x_r$ and let
   $a_i = x_iM'_\Delta$. Then a basis of the commutant
   $M'_\Delta$  is the set we denote by
   $\mathcal B'(M_\Delta)$ consisting of all elements
   $v$ of the form
   $$v = [x_i, x_j]^{{a_{j_1}}^{t_1}\ldots {a_{j_m}}^{t_m}}, \quad \{t_1,\ldots,t_m \}
     \subset \mathbb Z \setminus \{0\}$$
   such that the following conditions are satisfied:\\
   (1)
   $j \leq j_1 < j_2 \ldots < j_m \leq r,\,\,1 \leq j < i \leq r;$\\
   (2) the vertices
   $x_i$ and
   $x_j$ are in different connected components of
   $\Delta_v$ generated by all vertices of the set
   $\{x_i, x_j, x_{j_1},\ldots, x_{j_m}\};$\\
   (3)
   $x_i = \max (\Delta_{v,x_i})$, where
   $\Delta_{v,x_i}$ the connected component of the graph
   $\Delta_{v}$ containing
   $x_i$.
 \end{ttt}
 \begin{ccc}
   Let
   $\mathcal B'(M_\Delta)$ be linearly ordered. Then any element
   $g$ of the group
   $M_\Delta$ can be uniquely written in the form
   $$g =x_1^{\alpha_1}\ldots x_r^{\alpha_r}v_1^{\beta_1}\ldots v_m^{\beta_m},$$
   where
   $\alpha_i, \beta_j \in \mathbb{Z}$ and
   $v_1< \ldots < v_m, \,\,v_j \in \mathcal B'(M_\Delta)$.
 \end{ccc}
\subsubsection*{Centralizer dimensions}
 The notion of the centralizer dimension was introduced by
 A.\,Myasnikov and P.\,Shumyatsky
\cite{MS04} for comparison of universal theories of groups. Suppose that a sequence
 $$A_1 \subset A_2 \subset \ldots \subset A_n$$
 of subsets of a group
 $G$ is such that the chain of centralizers of these subsets
 $$C(A_1) > C(A_2) \ldots > C(A_n)$$
 is strictly descending. The \emph{centralizer dimension} of a group
 $G$ is the greatest
 $n$ for which such subsets
 $A_1,A_2,\dots, A_n$ of
 $G$ exist. A centralizer dimension is denoted by
 $Cdim(G)$. If the greatest
 $n$ does not exist  then we write
 $Cdim(G)=\infty$.

 It is known
\cite{MS04} that the centralizer dimension of a finitely generated metabelian group is
 finite. In papers
\cite{Tim17-1, Tim18}, properties of centralizer dimensions of partially commutative groups
 were studied and the exact value of
 $Cdim(M_\Delta)$, where
 $\Delta$ is a tree or a cycle, was found.

 \begin{ttt}\label{cdimcyc}
\cite{Tim17-1} Let
  $\Delta$ be a tree with at least 3 vertices. If
  $\Delta$ is a star then
  $Cdim (M_\Delta) = 3$. Otherwise,
  $Cdim(M_\Delta)=5$.
 \end{ttt}

 Let
 $g_1,\ldots, g_m \in M_\Delta$ be a finite system of elements and
 $\overline{g}_1,\ldots, \overline{g}_m$ the images of
 $g_1,\ldots, g_m$ in the free abelian group
 $\overline{M}_\Delta$ via the natural homomorphism
 $M_\Delta \to \overline{M}_\Delta$.   The rank of the system
 $g_1,\ldots, g_m$ is the rank of the subgroup generated by
 $\overline{g}_1,\ldots, \overline{g}_m$.

 Let
 $M_\Delta$ be a non-abelian group.
 Define a parameter
 $\alpha(M_{\Delta})$ for
 $M_\Delta$ as follows. Let
 $\mathcal{Z}(M_\Delta)=1$. Then put
 $\alpha(M_\Delta) = a$, where
 $a$ is the largest integer such that for any system of elements
 $g_1,\ldots, g_m \in M_\Delta$ of rank at least
 $a$ the centralizer
 $C(g_1,\ldots, g_m)$ is trivial. If the center of
 $M_\Delta$ is non-trivial then
 $M_\Delta/ \mathcal{Z}(M_\Delta)$ is a partially commutative group with no center. In this
 case, put
 $\alpha(M_\Delta) = \alpha(M_\Delta/ \mathcal Z(M_\Delta))$.

 Define a parameter
 $\beta$ for a group
 $M_\Delta$ as follows. Let
 $b$ be the least natural number such that for any distinct vertices
 $x_{i_1}, \ldots,x_{i_b}$ of
 $\Delta$ the intersection
 $$ \mathcal C(x_{i_1}) \cap \ldots \cap  \mathcal C(x_{i_b}) $$
 is trivial. Then put
 $\beta(M_\Delta) = b$.

 \begin{ttt}\label{cdim}
\cite{Tim18} Let
  $M_\Delta$ be a non-abelian group. Then
  $$Cdim(M_\Delta) \leq \alpha(M_\Gamma) + \beta(M_\Delta)+1.$$
 \end{ttt}

 \begin{ccc}
\cite{Tim18} For any group
  $M_\Delta$ the following equation holds:
  $$Cdim(M_\Delta) \leq 2n+1,$$
  where
  $n$ is quantity of vertices of
  $\Delta$.
 \end{ccc}

 The following theorem shows that the value of centralizer dimension is not bounded
 by a function of rank of a group and this value can be arbitrarily large even in the
 case of a two-generated metabelian group.

 \begin{ttt}
\cite{Tim18} For any
  $n \in \mathbb{N}$ there exists a two-generated untwisted metabelian group of
  centralizer dimension at least
  $n$.
 \end{ttt}

 By Theorem~%
\ref{cdimcyc}, centralizer dimensions of partially commutative metabelian groups defined by
 trees are bounded as well as centralizer dimensions of partially commutative groups defined
 by cycles.

 \begin{ttt}
\cite{Tim18} If
   $M_\Delta$ is a partially commutative group defined by a cycle of length at least
   $5$ then
   $Cdim(M_\Delta)=7$.
 \end{ttt}

 By analogy with
 $Cdim(G)$ the centralizer dimension in commutant
 $\mathcal{C}dim(G)$ is defined. To define
 $\mathcal{C}dim(G)$ centralizers in commutant
 $\mathcal{C}(A_i)=C(A_i)\cap G'$ are considered instead of centralizers
 $C(A_i)$. Centralizer dimensions in commutant can also be used for a comparison of
 universal theories.

 In
\cite{Tim17-1}, the centralizer dimensions
 $\mathcal{C}dim$ are defined for partially commutative groups defined by trees and cycles.

\subsubsection*{Inclusions, subgroups, retracts}
 For any variety
 $\mathfrak M$ and any graph
 $\Delta$ the following statement holds. Let
 $\Delta$ be a graph and
 $\Gamma$ its subgraph generated by a set of vertices
 $Y \subseteq V(\Delta)$. Then there exists a retraction of
 $G(\mathfrak M, \Delta)$ onto the group
 $G(\mathfrak{M}, \Gamma)$ such that this retraction is identical on
 $y$ and takes all other elements in
 $V(\Delta)$ to the unity.

 Let us describe a couple of embeddings of partially commutative metabelian groups into a
 group of matrices. They allied to the Magnus embedding. This embedding is very important
 in theory of soluble groups.  Recall the definition of the Magnus embedding for a free
 metabelian group
 $M_n$ for
 $n \geq 2$. Let
 $X=\{x_1, \ldots, x_n\}$ be a basis of
 $M_n$,
 $A_n$ a free abelian group with a basis
 $\{a_1,\ldots, a_n\}$,
 $B=\mathbb Z[A_n]$, and
 $F$ a free right
 $B$-module with a basis
 $\{f_1,\ldots, f_n\}$. Consider a matrix group
 $$W_n =
   \left(\begin{array} {ll}
    A_n & 0 \\
    F & 1
   \end{array}\right).$$
 The map
 $$\mu: x_i \mapsto  \left(
    \begin{array} {ll}
     a_i & 0 \\
     f_i & 1
    \end{array}\right), \text{ for } i = 1,\ldots,n,$$
 is extended up to the Magnus embedding
 $\mu$ of  the group
 $M_n$ to the group
 $W_n$.

 The embedding
 $\mu$ takes a commutator
 $[x_i, x_j]$ to the matrix
 $$\left(
   \begin{array} {ll}
     1 & 0 \\
     \tau_{ij} & 1
   \end{array}\right),\text{ where } \tau_{ij} = f_i(a_j-1) + f_j(1-a_i).$$
 Let
 $\Delta = (X; E)$ be a graph and
 $R_{\Delta}$ a normal subgroup generated by all commutators
 $[x_i, x_j]$ such that
 $(x_i, x_j) \in E$. Then
 $\mu$ maps
 $R_{\Delta}$ to the submodule
 $L$ of the module
 $F$ such that
 $L$ is generated by all
 $\tau_{ij}$ for which
 $(x_i, x_j) \in E$. Let
 $T = F/L$. The Magnus embedding
 $\mu$ of
 $M_n$ to
 $W_n$ induces an embedding
 $\mu_\Delta$ of the group
 $M_\Delta$ to the group of matrices
 \begin{equation} \label{magnus}
   W_\Delta = \left(
    \begin{array} {ll}
      A_n & 0 \\
      T & 1
    \end{array}\right).
 \end{equation}

 In
\cite{Tim18-1} the existence of one more embedding of a group
 $M_\Delta$ to a group of matrices was shown.

 \begin{ttt}
\cite{Tim18-1} Let
   $\Delta$ be a connected graph,
   $\{x_1,\ldots, x_n\}$ the set of vertices of
   $\Delta$ and a basis of the free metabelian group
   $M_n$,
   $a_i$ the image of
   $x_i$ under the natural homomorphism
   $M_n \to A_n = M_n/ M_n'$, and
   $\delta =(a_1 - 1)\cdot \ldots \cdot(a_n -1)$. Then the group
   $M_\Delta$ is embeddable to the group of matrices
   \begin{equation}\label{magnusfac}
     \overline{W}_\Delta = \left(
       \begin{array} {ll}
         A_n & 0 \\
         T/T\delta & 1
       \end{array}\right),
   \end{equation}
   where the module
   $T = F/L$ is defined above.
 \end{ttt}

 Groups
 $W_{\Delta}$ and
 $\overline{W}_{\Delta}$ are splittable. For this reason, they are preferable for a study
 of universal theories.  This will be discussed in Subsection
\ref{ssec12}.

 Let us present some theorems on subgroups. In
\cite{Tim10} it was shown that if a group
 $M_{\Delta}$ is nilpotent then it is abelian. The following theorem states even more.

 \begin{ttt}
\cite{Tim18-1} Any nilpotent subgroup of
   $M_{\Delta}$ is abelian.
 \end{ttt}

 The \emph{Fitting subgroup}
 $Fit(G)$ of a group
 $G$ is the product of all nilpotent normal subgroups of
 $G$.

 \begin{ttt}
\cite{Tim18-1} The Fitting subgroup of a group
  $M_\Delta$ is equal to the direct product of the center and the commutant of this group.
 \end{ttt}

 Let
 $G = (X\,|\,R,\mathfrak{A}^2)$ and
 $H = (Y\,|\,S,\mathfrak{A}^2)$ be represented in the variety of metabelian groups
 by generators and defining relations. If
 $X \cap Y = \varnothing$ then  the group
 $T = (X\sqcup Y\,|\, R \sqcup S,\mathfrak A^2)$ is called the metabelian product of
 $G$ and
 $H$. Let us denote the metabelian product of metabelian groups
 $G_1,\ldots, G_n$ by
 $M(G_1,\ldots, G_n)$.

 \begin{ttt}
\cite{Tim13} Any partially commutative metabelian partially commutative group is a subgroup
   of a direct product of some finite (possibly empty) set of free abelian groups
   $A_i$ and some finite (possibly empty) set of metabelian products
   $M_j=M(B^{(j)}_{1},\ldots, B^{(j)}_{{r_j}})$ of free abelian groups
   $B^{(j)}_{1},\ldots, B^{(j)}_{r_j}$.
 \end{ttt}

\subsubsection*{Automorphisms}
 A vertex
 $x$ of a graph
 $\Delta$ is called an end-point if its degree is equal to 1.

 An automorphism
 $\alpha$ of a group
 $G$ is called an
 $IA$-automorphism if this automorphism acts identically on the quotient group
 $\overline{G}$. The group of
 $IA$-automorphisms is denoted by
 $IAut(G).$

 \begin{ttt}\label{IAaut}
\cite{Tim20} Suppose that a graph
   $\Delta$ has no cycles. If an
   $IA$-au\-to\-mor\-phism
   $\alpha$ of the group
   $M_{\Delta}$ fixes all end-points and isolated vertices of
   $\Delta$ then $
   \alpha$ is the identical automorphism.
 \end{ttt}

 Both requirements are essential. If there is a cycle in
 $\Delta$ or
 $\alpha$ is not identical on the quotient group by the commutant then Theorem~%
\ref{IAaut} does not hold.

 Each automorphism
 $\alpha$ of a group
 $M_{\Delta}$ induces an automorphism
 $\overline{\alpha}$ of the free abelian group
 $\overline{M}_{\Delta}$. A group of induced automorphisms is called a group of quotient
 automorphisms of a group
 $\overline{M}_{\Delta}$ and is denoted by
 $\mathcal{F}_\Delta$. Clearly,
 $\mathcal F_\Gamma \cong Aut(M_\Gamma)/IAut(M_\Gamma)$.

 In
\cite{Tim20}, a description of a group of matrices
 $\mathcal{M}_\Delta$ is given. An automorphism
 $\overline{\alpha}$ of a group
 $\overline{M}_{\Delta}$ is a matrix automorphism if its matrix
 $[\overline{\alpha}]$ in some basis chosen by the graph
 $\Delta$ belongs to the group of matrices
 $\mathcal{M}_{\Delta}$. The following statement holds.

 \begin{ttt}
\cite{Tim20} Let
   $\Delta$ be a graph with no cycles. Then each quotient group automorphism of the group
   $M_\Delta$ can be written as a product of an automorphism of the graph
   $\Delta$ and a matrix automorphism.
 \end{ttt}

 Groups
 $G$ and
 $H$ are called \emph{commensurable} if there exist subgroups
 $G_1$ and
 $H_1$ of finite indices of the groups
 $G$ and
 $H$ respectively such that
 $G_1 \cong H_1$.

 Let a linear group
 $A$ be
 $\mathbb{Q}$-definable. This means that
 $A \leq GL(n,\mathbb{C})$ and its basic set is defined by a system of equations with
 coefficients in
 $\mathbb{Q}$. A subgroup
 $B \leq A\, \cap\, GL(n,\mathbb Q) = A_{\mathbb Q}$ is called an arithmetic group or an
 arithmetic subgroup of
 $A$ if it is commensurable with
 $A_{\mathbb Z} = A\cap GL(n,\mathbb Z)$.

 \begin{ccc}
\cite{Tim20}. Let
  $\Gamma$ be a graph with on cycles. Then the group of quotient automorphisms of
  $\overline{M}_{\Delta}$ is arithmetic.
 \end{ccc}

 Let us give some more information on automorphisms of partially commutative metabelian
 groups (see
\cite{Tim15} for details). Let
 $\Delta=(X; E)$, be a graph with the set of vertices
 $X=\{x_1, \ldots, x_n\}$. In
\cite{Lau95}, Laurence defined four sets of automorphisms generating the group
 $Aut(F_{\Delta})$.\\
 (1) The set of graph automorphisms, namely the set of elements in
 $Aut(F_{\Delta})$ such that these elements are induced by automorphisms
 $\pi: \Delta \to \Delta$ of the graph
 $\Delta$.\\
 (2) The set of inverting automorphisms
 $\alpha \in Aut(F_\Delta)$. These are automorphisms taking one of the vertices
 $x_i \in X $ to
 $x_i^{-1}$ and fixing all other vertices.\\
 (3) Consider two distinct vertices
 $x_i, x_j$, such that
 $(x_j, x) \in E$ implies
 $(x_i, x) \in E$ for any
 $x \in X$. The third set consists of transvections taking
 $x_j$ to
 $x_jx_i^{\pm 1}$ or to
 $x_i^{\pm1}x_j$ and fixing all other vertices.\\
 (4) The fourth set consist of locally interior automorphisms defined as follows. Let
 $x_i \in X$. Consider the subgraph
 $\Gamma$ obtained by deleting
 $x_i$, all vertices adjacent to
 $x_i$, and all edges incident to deleted vertices. Let
 $\Lambda$ be a union of some connected components of
 $\Gamma$. Then define
 $\beta \in Aut(F_{\Delta})$ setting
 $\beta(x_j) = x_i^{-1}x_j x_i$ for
 $x_j\in \Lambda$ and
 $\beta(x_j) = x_j$ for
 $x_j \notin \Lambda.$

 It follows from
\cite{Ush01} that if
 $\Delta$ has at least three connected components and each of these components is a complete
 graph then the group
 $Aut(M_{\Delta})$ is not generated by automorphisms induced by the Laurence automorphisms.

 In
\cite{Tim15}, a stronger result was obtained. Namely, if
 $\Delta$ is a connected graph or even a tree then the group
 $Aut(M_\Delta)$ can contain automorphisms not induces by automorphisms of partially
 commutative group
 $F_{\Delta}$. In the same paper,  a monoid
 $\mathcal{P}$ of matrices over a ring
 $\mathbb{Z}[a_1^{\pm 1},\ldots, a_n^{\pm 1}]$ of integer Laurent polynomials and a
 congruence
 $\approx$ on this monoid are defined in such a way that the group of automorphisms acting identically
 on the quotient group
 $\overline{M}_{\Delta}$ was isomorphic to the quotient monoid of
 $\mathcal{P}$ by
 $\approx$.

 The structure of the group of automorphisms of partially commutative class two nilpotent
 group was studied in
\cite{RT10}.

\subsubsection*{Direct decompositions}
 A group
 $G$ is \emph{decomposable into a direct product} if
 $G=A \times B$ for some groups
 $A \neq 1$ and
 $B \neq 1$. In
\cite{RT20}, the question on existence of direct decomposition for partially commutative
 groups in varieties containing
 $\mathfrak N_2$ was studied.  Two theorems were proved.

 \begin{ttt}
\cite{RT20} Let
  $\mathfrak{M}$ be a variety of groups such that this variety contains
  $\mathfrak N_2$. Suppose that a group
  $G = G(\mathfrak{M},\Gamma)$ decomposes into a direct product
  $H \times A$, where
  $A$ is an abelian group. Then there exists a subgraph
  $\Delta$ of the graph
  $\Gamma$ such that the set of vertices of
  $\Delta$ contains
  $X\backslash X^{\perp}$ and
  $H\cong G(\mathfrak{M},\Delta)$.
 \end{ttt}

 \begin{ttt}
\cite{RT20}. Let
  $\mathfrak{M}$ be a variety of soluble groups such that this variety contains
  $\mathfrak{N}_2$. If a graph
  $\Delta$ is not connected then the group
  $G(\mathfrak M,\Delta)$ is not decomposable into a direct product.
 \end{ttt}

\subsection{Logical properties of partially commutative metabelian groups} \label{ssec12}

\subsubsection*{Elementary equivalence and isomorphism}
 In
\cite{KMNR80}, it was shown that two partially commutative associative algebras are
 isomorphic if and only if so are their defining graphs. Using this result, C.\,Droms
\cite{Dro87} proved an analogous one for partially commutative groups in the variety
 of all groups. He established the following fact.

 \begin{ttt}
\cite{Dro87} If
  $\mathfrak{M}$ contains the variety
  $\mathfrak{N}_2$ then
  $G(\mathfrak{M},\Gamma)$ and
  $G(\mathfrak{M},\Delta)$ are isomorphic if and only if so are the graphs
  $\Gamma$ and
  $\Delta$.
 \end{ttt}

 The following theorem provides a criterium of coincidence of elementary theories
 of a partially commutative group in a nilpotent variety containing
 $\mathfrak{N}_2$ and an arbitrary group.

 \begin{ttt}\label{eleqpcgg}
\cite{RT20} Let a variety of nilpotent groups
  $\mathfrak{M}$ contain
  $\mathfrak{N}_2$. Then if a finitely generated group
  $H$ has the same elementary theory as
  $G(\mathfrak M,\Delta)$  then  $G(\mathfrak M,\Delta) \cong H$.
 \end{ttt}

 Theorem~%
\ref{eleqpcgg} implies the following result.

 \begin{ccc}
\cite{RT20} Let
  $G=G(\mathfrak M,\Gamma)$ and
  $H=G(\mathfrak M, \Delta)$ be groups in a variety
  $\mathfrak{M}$ of nilpotent groups such that
  $\mathfrak{M}$ contains
  $\mathfrak{N}_2$. Then the elementary theories of the groups
  $G$ and
  $H$ coincide if and only if
  $\Gamma \cong \Delta.$
 \end{ccc}

 Let
 $Th(G)$ be an elementary theory. An elementary theory
 $Th(G)$ is called \emph{soluble} if there is an effective procedure checking for any
 sentence
 $\Phi$ if this sentence belongs to
 $Th(G)$.

 In
\cite{Nos12}, G.\,A.\,Noskov proved that the elementary theory of an almost soluble
 group is soluble if and only if the group is almost abelian. So, if a variety
 $\mathfrak{M}$ is soluble and
 $\mathfrak{N}_2 \subseteq \mathfrak{M}$ then the elementary theory of a group
 $G(\mathfrak M, \Delta)$ is not soluble.

\subsubsection*{Universal theories}

 One of the reasons of making researches of centralizer dimensions
 $Cdim$ and
 $\mathcal{C}dim$ is a coincidence of universal theories of groups implies an equality
 of their centralizer dimensions.

 In Subsection
\ref{ssec11}, inclusions
 \eqref{magnus} and
\eqref{magnusfac} of a group
 $M_{\Delta}$ to the groups of matrices
 $W_{\Delta}$ and
 $\overline{W}_{\Delta}$ were defined. In
\cite{GT09}, it was shown that the universal theories of these groups of matrices are soluble. But the universal
 theory of
 $M_{\Delta}$ coincides with no universal theories of groups of matrices. This result was
 obtained in
\cite{Tim18}. Therefore we can only say that a group
 $M_{\Delta}$ is embeddable into a metabelian group with a soluble universal theory.
 The problem on the solubility of the universal theory of a group
 $M_{\Delta}$ has not been solved yet. It is included into the Kourovka Notebook
\cite{Kourovka}. It is known that the universal theory of a free abelian groups is
 soluble. In
\cite{Cha95, Cha97}, O.\,Chapius proved that the universal theory of a free metabelian group
 is also soluble. Obviously, groups with soluble universal theories can be obtained from free abelian and free metabelian groups by using the direct product of groups. For instance, if
 $\Delta_4$ is the
 $4$-cycle then the universal theory of the group
 $G_{\Delta_4}$ is soluble, since it is isomorphic to the direct product of two free
 $2$-generated metabelian groups. One can find a non-trivial example of a partially
 commutative metabelian group having a soluble universal theory. So, in
\cite{Tim17-1}, it was shown that if
 $\Gamma_4$ is the linear graph on four vertices then the universal theory of partially
 commutative metabelian group
 $M_{\Gamma_4}$ is soluble. The proof follows from the coincidence of universal theories of groups $M_{\Delta_4}$ and
 $M_{\Gamma_4}$,

 In
\cite{BT17}, the problem on the universal equivalence of partially commutative metabelian
 groups with acyclic defining graphs was considered. The following theorem was proved.

 \begin{ttt}\label{abelguni}
\cite{BT17} Suppose that the graph
  $\Delta^*$ is obtained from a graph
  $\Delta$ by deleting all end-points and the edges incident to the end-points. Let
  $\mathfrak{A}_p$ be the variety of abelian groups of exponent
  $p$, where
  $p$ is a prime number or
  $0$. If
  $\Delta$ and
  $\Gamma$ are graphs with no cycles such that each connected component of these graphs
  has at least three vertices then the groups
  $G(\mathfrak{A}_p \mathfrak{A}, \Gamma)$ and
  $G(\mathfrak{A}_p \mathfrak{A}, \Delta)$ are universally equivalent if and only if
  $\Gamma^* \cong \Delta^*$.
 \end{ttt}

 If any connected component of
 $\Delta$ of
 $\Gamma$ contains less then tree vertices then the corresponding connected component in
 $\Delta^*$ (correspondingly in
 $\Gamma^*$) is empty and the statement of Theorem~%
\ref{abelguni} does not hold.

 To prove Theorem~%
\ref{abelguni}, generalizations many algebraic properties of a group
 $M_{\Delta}$ to groups
 $G(\mathfrak{A}_p \mathfrak{A},\Delta)$ were used. These generalizations were obtained in
\cite{BT17}.

 A study of partially commutative metabelian group were continued in
\cite{GT11-1}. In this paper, an equivalence relation on the set of vertices of a graph
 $\Delta$ is defined. Then, an adjacency relation is determined on the set of equivalence classes. So,
 a new graph
 $\Delta^\star$, appears.  This graph is called the \emph{compression} of the initial one
 and it is usually simpler then
 $\Delta$. Let us give a strict definition of the compression of a graph. We say that two
 vertices
 $x$ and
 $y$ of a graph
 $\Gamma$ are equivalent and write
 $x \sim y$ if
 $x^\bot \cup \{x\} = y^\bot \cup \{y\}$
 ($x^\bot$ is defined in~%
\eqref{orthogonal}). Note that equivalent vertices are always adjacent. Then the compression
 of
 $\Delta$ is a quotient graph
 $\Delta^\star =\Delta/\sim$.

 An end-point
 $z$ in
\cite{GT12} is called \emph{bad} if there exists a vertex
 $x$ adjacent to
 $z$ and at least to two vertices
 $y$ and
 $v$, where
 $y$ is also an end-point.

 Denote by
 $\Delta'$ a graph obtained by deleting one-by-one all bad vertices
 and the edges incident to them.

 \begin{ttt}\label{starmaltese}
\cite{GT11-1,GT12} For any graph
  $\Delta$ the universal theories of the groups
  $M_{\Delta}$ and
  $M_{{\Delta}^\star}$ coincide as well as the universal theories of
  $M_{\Delta}$ and
  $M_{\Delta'}$.
 \end{ttt}

 The paper
\cite{GT11-1} gives an example of  graph
 $\Delta$ such that this graph is not a tree while its compression
 $\Delta^\star$ is. For this reason, the universal theory of a group defined by a tree can
 coincide with the universal theory of a group defined by a graph with cycles. Theorem~%
\ref{starmaltese} implies that the condition of acyclicity of a defining graph in Theorem~%
 \ref{abelguni} is essential even for partially commutative groups in the variety
 $\mathfrak A^2.$

 It follows from Theorem~%
\ref{starmaltese} that the universal theory of metabelian product of free abelian groups
 coincides with one of a free metabelian group.

 For partially commutative metabelian groups defined by cycles the following theorem holds.

 \begin{ttt}
\cite{GT11-1}. If
  $n,m \geq 3$ then groups
  $M_{\Delta_n}$ and
  $M_{\Delta_m}$ defined by cycles of lengths
  $n$ and
  $m$ respectively are universally equivalentif and only if
  $n=m$.
 \end{ttt}

 Let
 $\Delta$ be a graph with the set of vertices
 $X= \{x_1,\dots,x_n\}$. Denote by
 $\varphi (\Delta)$ the following sentence.

 $$\exists v_1 \ldots v_n (\bigwedge_{(x_i, x_j) \in \Delta}[v_i,v_j]=
    1 \wedge \bigwedge_{(x_i, x_j) \notin\Delta}[v_i, v_j] \neq 1 \wedge
    \bigwedge_{i\ne j} v_i \neq v_j\wedge \bigwedge_{i=1,n}v_i \neq 1).$$

 V.\,N.\,Remeslennikov formulated the following conjecture. \\[1ex]
 {\it Let
 $\mathfrak{M}$ be a variety of groups. If the universal theories of the groups
 $F(\mathfrak M, \Gamma_1)$ and
 $F(\mathfrak M, \Gamma_2)$ are distinct then there exist a graph
 $\Delta$ such that the sentence
 $\varphi (\Delta)$ is true on one of these groups and is false on the other one.}\\[-1ex]

 In
\cite{MT07}, the affirmative solution of this conjecture was obtained for partially
 commutative nilpotent
 $R$-groups of class 2, where
 $R$ is a binomial ring. Let
 $G$ and
 $H$ be two partially commutative nilpotent
 $R$-groups of class 2 and
 $\Delta$ and
 $\Gamma$ defining graphs of
 $G$ and
 $H$ respectively. It turns out that if
 $G$ and
 $H$ are not universally equivalent then their universal theories differ in
 $\varphi(\Gamma)$ or
 $\varphi(\Delta)$.

 However, for partially commutative metabelian groups the analogous result
 does not hold. A counterexample was obtained in
\cite{Tim17-1}. It is not known if the conjecture holds for the variety of
 metabelian groups. Nevertheless, if only formulas of the form
 $\varphi (\Delta)$, where
 $\Delta$ is a tree, are considered then the corresponding result is not true. In
\cite{Tim17-1}, the second author of this survey has found two groups
 $M_{\Gamma_1}$ and
 $M_{\Gamma_2}$ such that these groups have distinct universal theories while for
 any tree
 $T$ the corresponding formula
 $\varphi (T)$ is true on one of these groups if and only if it is true on the
 other one.

\subsubsection*{Quasi-varieties}
 A sentence of the type
 $$\forall z_1\ldots, z_m ((w_1(z_1,\ldots, z_m)=1 \wedge \ldots
   \wedge w_r(z_1,\ldots,z_m)=1) \longrightarrow w(z_1,\ldots,z_m)=1),$$
 where
 $w$ and
 $w_i$ are group words is called a \emph{quasi-identity}. A class of groups satisfying
 a collection of quasi-identities is called a \emph{quasi-variety}. A non-empty class of
 groups form a quasi-variety if and only if this class is closed with respect to taking
 subgroups, cartesian products, and ultra-products. Denote by
 $qvar(G)$ the quasi-variety generated by a group
 $G$. A class of groups closed with respect to taking subgroups and cartesian products is
 called a pre-variety.

 In
\cite{Tim13} it was shown that there exist free partially commutative groups
 $F_{\Delta_i}$  such that
 $$qvar (F_{\Delta_1}) \subset qvar (F_{\Delta_2})\subset \ldots \subset qvar (F_{\Delta_n})
   \subset \ldots ,$$
 and all inclusions in this infinite chain are strict. This is not so for partially
 commutative  metabelian groups. For them, the following theorem holds.

 \begin{ttt}
\cite{Tim13} Any two non-abelian partially commutative metabelian groups generate equal
   quasi-varieties
 \end{ttt}

 The same result takes place for pre-varieties.

 \begin{ttt}
\cite{Tim13} Any two non-abelian partially commutative metabelian groups generate equal
  pre-varieties.
 \end{ttt}

 The \emph{positive universal theory} of a group
 $G$ is the set of all sentences
 $\Phi$ of the form
 $$\Phi = \forall z_1\ldots z_m \left(\bigvee_{i \in I}\bigwedge_{j_i \in J_i}w_{j_i}
   (z_1,\ldots,z_m)=1\right),$$
 such that these sentences are true on
 $G$, where
 $w_{j_i}(z_1,\ldots,z_m)$ are group words.

 Let us denote the positive universal theory of a group
 $G$ by
 $Th^+_{\forall}(G).$

 The following theorem shows that not only quasi-varieties and pre-varieties generated by
 partially commutative metabelian groups coincide but also positive universal theories of
 such groups do.

 \begin{ttt}
\cite{Tim13} Let
  $M_{\Gamma}$ and
  $M_{\Delta}$ be non-abelian groups. Then
  $Th^+_{\forall}(M_{\Gamma})=Th^+_{\forall}(M_{\Delta})$.
 \end{ttt}

\subsubsection*{Equations}
 In Subsection~%
\ref{ssec11}, the inclusion of a group
 $M_{\Delta}$ into the corresponding group of matrices
 $W_{\Delta}$ was defined. In
\cite{Tim18-2}, it was shown that, in general, the universal theories of groups
 $M_{\Delta}$ and
 $W_{\Delta}$ are distinct. This result was obtained by comparing centralizer
 dimensions of these groups. However, the groups
 $M_{\Delta}$ and
 $W_{\Delta}$ have some common properties allied to their universal theories. Namely, the
 following theorem was proven.

 \begin{ttt}\label{equat}
\cite{Tim18-2} An equation
   $$g_1x^{m_1} \ldots g_lx^{m_l} = 1,\,\,g_i \in M_\Delta,$$
   is solvable in
   $M_{\Delta}$ if and only if it is so in
   $W_{\Delta}$.
 \end{ttt}

 The analogue of Theorem~%
\ref{equat} does not hold for equations of two unknowns. Moreover,
 this analogue does not hold even for a totally disconnected graph
 $\Delta.$

\subsection{Partially equivalent metabelian nilpotent groups} \label{ssec13}
\subsubsection*{Mal'cev basis}
 Let
 $\mathfrak{N}_{2, c}$ be the intersection of the variety of metabelian groups
 $\mathfrak{A}^2$ with the variety of nilpotent groups
 $\mathfrak{N}_c$.

 We introduce some notation from
\cite{Tim11}. In this paper, a basis of a group
 $M_{c,\Delta}=G(\mathfrak N_{2, c}, \Delta)$ is constructed. Let
 $G$ be a  finitely generated nilpotent torsion-free group. As it is known,
 $G$ has a central series
 $$G=G_1> G_2>\ldots>G_{s+1} =1$$
 with infinite cyclic quotient groups. Let us choose elements
 $a_1,\ldots, a_s \in G$ such that
 $G_i = \langle a_i,\,\,G_{i+1}\rangle$. An ordered system of elements
 $(a_1,\ldots, a_s)$ is called a \emph{Mal'cev basis} of a group
 $G$. Each element
 $g \in G$ can be written in the form
 $$g = a_1^{t_1}\ldots a_s^{t_s},\quad t_i \in \mathbb {Z}$$
 uniquely.

 Let
 $\Delta=(X;E)$ be a graph with the set of vertices
 $X =\{x_1,\ldots, x_n\}$ and
 $v(x_{i_1},\ldots, x_{i_m})$ a representation of an element
 $v \in M_{c,\Delta}$ via generators in
 $X$, where vertices
 $x_{i_1}, \ldots, x_{i_m}$ occur in this representation. Then set
 $\sigma(v) = \{ x_{i_1}, \ldots, x_{i_m}\}$. Note that,
 $\sigma(v)$ depends not only on
 $v$ but also on a specific representation via generators of the group. Denote by
 $\Delta_v$ the subgraph of
 $\Delta$ generated by the set
 $\sigma(v)$. The connected component of the graph
 $\Delta_v$ containing a vertex
 $x \in \sigma(v)$ is denoted by
 $\Delta_{v,x}$. Let us order the set
 $X$ as
 $x_1 < x_2 <\ldots < x_r$. Denote the greatest vertex in a connected component
 $\Delta_{v,x}$ by
 $\max(\Delta_{v,x})$. Define a commutator
 $c_m = [y_1, y_2,\ldots, y_m]$ by induction:
 $c_2 = [y_1, y_2]$,
 $c_m = [c_{m-1}, y_m]$. Let
 $\mathcal{B}'( M_{c,\Delta})$ be the set of commutators of the form
 $$v = [ x_{j_1}, x_{j_2}, \ldots, x_{j_m}], \quad 2 \leq m \leq c$$
 in a group
 $M_{c,\Delta}$ such that the following conditions are satisfied:\\
 (1)
 $1 \leq j_2 \leq j_3 \leq \ldots \leq j_m \leq r,\,\,j_2 < j_1 \leq r$;\\
 (2) the vertices
 $x_{j_1}$ and
 $x_{j_2}$ are in different connected components of the graph
 $\Delta_v$;\\
 (3)
 $x_{j_1} = \max(\Delta_{v, x_{j_1}})$.

 \begin{ttt}\label{malbas}
\cite{Tim11} The set of elements
  $\mathcal{B}(M_{c,\Delta}) = X\sqcup\mathcal{B}'(M_{c,\Delta})$ is a Mal'cev
  basis of a group
  $M_{c,\Delta}$.
 \end{ttt}

 The canonical representation from Theorem~%
\ref{malbas} is used in
 \cite{GT12} for study of algebraic properties and the universal theory of a group
 $G=M_{c,\Delta}$. Let us present the main results of this paper. The following theorem
 is similar to Theorem~%
\ref{annihilators} on annihilators of partially commutative metabelian groups and it uses
 the ideals
 $\mathcal A_{ij}$ defined in Theorem
\ref{annihilators}. Let
 $\triangle$ be a fundamental ideal of a ring
 $\mathbb Z[\overline{G}]$, i.e. the kernel of the natural homomorphism
 $\mathbb Z[\overline{G}] \to \mathbb Z$.

 \begin{ttt}
\cite{GT12} Let
  $x_i$ and
  $x_j$ be two non-adjacent vertices of a graph
  $\Gamma$. Then the annihilator of the commutator
  $[x_i,x_j]$ in a group
  $M_{c,\Delta}$ is equal to
  $\mathcal{A}_{i,j} + \triangle ^ {c-1}$.
 \end{ttt}

 Let us present a theorem on the centralizers for elements of a group
 $M_{c,\Delta}$. Denote by
 $C(g)$ the centralizer of
 $g$ and by
 $\mathcal{C}(g)$ the centralizer of
 $g$ in commutant, namely the set
 $$\mathcal{C}(g)= C(g) \cap  M'_{c,\Delta}.$$
 As usual, let
 $\gamma_m(G)$ denote the
 $m$th element of the lower central series of
 $G$.

 \begin{ttt}
\cite{GT12} Let
  $X=\{x_1,\ldots, x_n\}$ be the set of vertices of a graph
  $\Delta$ and
  $G = M_{c,\Delta}$. Then the following conditions hold.\\
  1) If
  $x_n$ in an isolated vertex then
  $$C(x_n) = \langle x_n\rangle \times \gamma_c(G).$$
  2) If
  $x_n$ is adjacent to only one vertex (ex., to
  $x_{n-1}$) then
  $$C(x_n) = \langle x_{n-1}\rangle \times \langle x_{n}\rangle \times \gamma_c(G).$$
  3) If
  $x_n$ is adjacent to vertices
  $x_{r+1},\ldots,x_{n-1}$, where
  $r\leq n-3$, then
  $C(x_n)$ consists of all elements of the form
  $$\prod_{p=r+1}^n {x_p}^{l_p} \cdot
    \prod_{r+1 \leq i < j \leq n-1}{[x_i,\, x_j]} ^{f_{i,j}} \cdot \gamma_c(G),$$
  where
  $l_p \in \mathbb Z,\,\, f_{i,j} \in \mathbb Z[\overline{G}]$.
 \end{ttt}

 The following theorem can be used to find the centralizers in commutant for any elements
 of a group
 $M_{c,\Delta}.$

 \begin{ttt}
\cite{GT12} Let
  $\Delta = (X; E)$ be a graph with the set of vertices
  $X=\{x_1,\ldots,x_n\}$ and
  $G=M_{c,\Delta}$. Then
  $$\mathcal C(gx_{i_1}^{l_1}\ldots x_{i_t}^{l_t}) = \bigcap_{j=1}^t \mathcal C(x_{i_j})$$
  for any integers
  $1 \leq i_1 < \ldots <i_t \leq n$, any non-zero integers
  $l_1,\ldots, l_t$, and any element
  $g$ in the commutant
  $G'$.
 \end{ttt}

 To study the universal theory of a group
 $M_{c,\Delta}$ the following theorem is useful.

 \begin{ttt}\label{intcentrnilpmetgr}
\cite{GT12} Let
  $\Delta = (X; E)$  be a tree,
  $u,v$ two distinct vertices of
  $X$. Then in the group
  $G=M_{c,\Delta}$ the following identity holds:
  $$\mathcal{C} (u)\cap \mathcal{C}(v) = \gamma_c(G).$$
 \end{ttt}

 Theorem~%
\ref{intcentrnilpmetgr} implies a description of the center of a partially commutative
 nilpotent metabelian group.

 \begin{ttt}
\cite{GT12} Let
  $\Delta = (X; E)$. Then the center of
  $G = M_{c,\Delta}$ is the direct product of the group
  $\gamma_c(G)$ and the cyclic groups generated by vertices
  $x_{i} \in X^\bot$.
 \end{ttt}

 Let us move on to the results on the universal theory of a group
 $M_{c,\Delta}$. The definition of a graph
 $\Delta'$ is given in Subsec.~%
\ref{ssec12} before Theorem~%
 \ref{starmaltese}, a graph
  $\Delta^*$ is defined in Theorem~%
\ref{abelguni}.

 \begin{ttt}
\cite{GT12} Let
   $\Delta$ be a graph. The groups
   $M_{c,\Delta}$ and
   $M_{c,\Delta'}$ have the same universal theories.
 \end{ttt}

 \begin{ttt}\label{metabnilgpuni}
\cite{GT12}. Let
  $\Gamma_1$ and
  $\Gamma_2$ be trees. The groups
  $M_{c,\Gamma_1}$ and
  $M_{c,\Gamma_2}$ have the same universal theories if and only if the graphs
  $\Gamma_1^*$ and
  $\Gamma_2^*$ are isomorphic.
 \end{ttt}

 Theorem~%
\ref{metabnilgpuni} is an analogue of Theorem~%
 \ref{abelguni}.

\subsection{Partially commutative metabelian
 pro-$p$-groups}\label{ssec14}

 In
\cite{AT19}, centralizers of elements and annihilators of commutators in
 partially commutative metabelian
 pro-$p$-groups were studied. The results obtained for partially commutative metabelian
 pro-$p$-groups are similar to those for partially commutative metabelian abstract groups in
\cite{GT09}. In this subsection, we are talking about
 pro-$p$-groups. So, by a subgroup, a homomorphism, a generating set we mean a closed
 subgroup, a continuous homomorphism, a generated set in topological sense respectively. Denote by
 $P$ a free metabelian
 pro-$p$-group and by
 $P_{\Delta}$ the partially commutative metabelian
 pro-$p$-group  defined by a graph
 $\Delta = (X;E)$. Let
 $X = \{x_1,\ldots, x_n\}$. A quotient group of
 $P_{\Delta}$ by its commutant
 $P'_{\Delta}$ is a free abelian
 pro-$p$-group
 $A$ with a basis
 $\{a_1,\ldots, a_n\}$, where
 $a_i$ is an image of
 $x_i$ via the natural homomorphism
 $P_{\Delta} \to P_{\Delta}/P'_{\Delta}$. This group is isomorphic to the direct sum of
 $n$ copies of additive group of the ring of integer
 $p$-adic numbers
 $\mathbb{Z_p}$. The action of
 $P_{\Delta}$ on
 $P'_{\Delta}$ by conjugations
 $$x\rightarrow x^g=g^{-1}xg$$
 defines a structure of a right module on
 $P'_{\Delta}$ over the augmented group algebra
 $\mathbb{Z}_p[[A]]$. This algebra is identified with the power series algebra
 $\mathbb Z_p[[y_1,\ldots, y_n]]$, where
 $y_i=a_i-1$. Similarly,
 $P'$ is a module over the algebra
 $\mathbb Z_p[[y_1,\ldots, y_n]]$. For this reason, any element
 $f \in P$ can be written in the form
 $$f = x_1^{l_1} \ldots x_n^{l_n} \prod_{1 \leq i <  j \leq n}
   [x_i,x_j]^{\alpha_{ij}},$$
 where
 $l_i \in \mathbb Z_p,\,\,\,\alpha_{ij} \in \mathbb Z_p[[y_1,\ldots, y_n]].$

 For a graph
 $\Delta$ and any its vertices
 $x_i$ and
 $x_j$, let us define the ideal
 $\mathcal A_{i,j}$ of the algebra
 $\mathbb Z_p[[y_1,\ldots, y_n]]$ as it was made for partially commutative metabelian
 group and for partially commutative metabelian nilpotent group. Namely, if the vertices
 $x_i$ and
 $x_j$ lie in different connected components of
 $\Delta$ then set
 $\mathcal{A}_{i,j} = 0$. If the vertices
 $x_i$ and
 $x_j$  lie in the same connected component then consider each path with no returns
 $(x_i, x_{i_1},\ldots,x_{i_r}, x_j)$ between these vertices. To each such path assign the
 product
 $y_{i_1} \ldots y_{i_r}$ if length of the path is greater than 1 and 1 otherwise. By
 definition, the ideal
 $\mathcal A_{i,j}$ is generated by all such elements. In particular, if
 $x_i, x_j$ are adjacent then
 $\mathcal{A}_{i,j}$ contains 1.  So, this ideal coincides with the entire algebra
 $\mathbb Z_p[[y_1,\ldots,y_n]].$

 Let us formulate the main results of paper
\cite{AT19}.

 \begin{ttt}
\cite{AT19} Let
  $x_1,\ldots, x_n$, where
  $n\geq 2$, be vertices of the defining graph
  $\Delta$ of a partially commutative metabelian
  pro-$p$-group
  $P_{\Delta}$. Then for
  $i \neq j$, the annihilator of the commutator
  $[x_i, x_j]$ in the algebra
  $\mathbb Z_p[[y_1,\ldots,y_n]]$ coincides with the ideal
  $\mathcal A_{i,j}$.
 \end{ttt}

 \begin{ttt}
\cite{AT19} Let
   $x_1, \ldots, x_n$, where
   $n\geq 2$, be vertices of the defining graph
   $\Delta$ of a partially commutative metabelian
   pro-$p$-group
   $P_{\Delta}$ and
   $x_2,\ldots, x_m$ all vertices adjacent to
   $x_1$. An element
   $g \in P_{\Delta}$ lies in the centralizer of
   $x_1$ if and only if it can be written in the form
   $$g = x_1^{l_1}\ldots x_m^{l_m} \prod_{2\leq i <j \leq m}[x_i, x_j]^{\gamma_{i,j}},$$
   where
   $l_i \in \mathbb Z_p,\,\,\gamma_{i,j} \in \mathbb Z_p[[y_1,\ldots,y_n]]$.
 \end{ttt}

 The following theorem has not been published yet. It is analogous to Theorem~%
\ref{bas_comm} and gives a description of a basis for the commutant of a partially
 commutative metabelian
 pro-$p$-group.

 \begin{ttt}
  Let the set
  $\{x_1,\ldots, x_n\}$ of vertices of a graph
  $\Delta$ be  ordered.  Then a basis
  $\mathcal{B}(P'_\Delta)$ of the commutant
  $P'_\Delta$ over
  $\mathbb{Z}_p$ is the set of all elements
  $w$ of the form
  $$w = [x_i, x_j]^{y_{j_1}^{s_1}\ldots y_{j_m}^{s_m}}, \{s_1,\ldots
     , s_m\} \subset \mathbb{N},$$
  such that the following conditions are satisfied:\\
  1) $x_j \leq x_{j_1}  < \ldots < x_{j_m}, \,\, x_j < x_i$;\\
  2) the vertices
  $x_i, x_j$ are in different connected components of graph
  $\Delta_w$ generated by all vertices of the set
  $\{x_i,x_j, x_{j_1},\ldots, x_{j_m}\}$;\\
  3) $x_i = \max \{\Delta_{w,x_i}\}$, where
  $\Delta_{w,x_i}$ the connected component of the graph
  $\Delta_w$ containing
  $x_i$.
 \end{ttt}


  \section{Partially commutative Lie algebras} \label{sec2}
  Researches of partially commutative Lie algebras began only about 30 years ago and this 
  algebras are studied not so heavily as partially commutative groups.

  Since partially commutative groups and Lie
  algebras are rather similar objects, many results for groups have analogues for Lie
  algebras. Moreover, some methods for studying Lie algebras come from ones for
  studying groups. Nevertheless, there are specific methods for researching Lie algebra.

  Let
  $\Gamma=(A;E)$ be a graph with a (finite or infinite) set  of vertices
  $A=\{a_1,a_2,\dots\}$ and a set of edges
  $E$. If
  $a_i$ and
  $a_j$ are adjacent in
  $\Gamma$ then we write
  $a_i\leftrightarrow a_j$. Similarly, if
  $B \subseteq A$ and
  $a_i\leftrightarrow a_j$ for any
  $a_j \in B$ then we write
  $a_i \leftrightarrow B$. Finally,  let
  $B,C\subseteq A$. Then
  $B \leftrightarrow C$ means
  $a_i \leftrightarrow a_j$  for any
  $a_i\in B$ and any
  $a_j \in C$.

  Given
  $B\subseteq A$ denote by
  $\Gamma(B)$ the full subgraph of
  $\Gamma$ generated by the set of vertices
  $B$, namely put
  $\Gamma(B)= (B;E\cap B^2)$. For a subgraph
  $\Delta$ of a graph
  $\Gamma=\langle A;E\rangle$ denote by
  $A(\Delta)$ the set of vertices of
  $\Delta$.

  Let
  $R$ be a unital commutative ring and
  $\Gamma=(A;E)$ be an undirected graph without loops. The definition of a partially
  commutative Lie
  $R$-algebra in a variety
  $\mathfrak{M}$ can be written as follows.
  $$L_R(\mathfrak{M}, \Delta)= \langle A \,|\, [a_i,a_j]=0
    \text{ if } a_i\leftrightarrow a_j; \mathfrak{M} \rangle.$$
  Indeed, by
\eqref{pcs}
  \begin{equation} \label{pcidlialg}
    [a_i,a_j]=[a_j,a_i]
  \end{equation}
  for any  \protect
  $i$ and
  $j$ such that
  $a_i \leftrightarrow a_j$. On the other hand,
  $[f,g]=-[g,f]$ for any elements
  $f$ and
  $g$ of any Lie
  $R$-algebra. In particular,
  \begin{equation} \label{ac}
    [a_i,a_j]=-[a_j,a_i].
  \end{equation}
  Combining
\eqref{pcidlialg} and
 \eqref{ac} we obtain
  $2[a_i,a_j]=0$ if
  $a_i \leftrightarrow a_j$. Therefore, characteristic of
  $R$ is equal to
  $2$ or
  $[a_i,a_j]=0$. In the former case, the corresponding Lie
  $R$-algebra is commutative, so the notion of partial commutativity makes no sense.
  So, in this section we assume that characteristic of a basic ring (or a field)
  is not equal to 2.

  In this section, we talk mainly about results in partially commutative and
  partially commutative metabelian Lie 
  $R$-algebras. Some results concern partially commutative  nilpotent Lie
  $R$-algebras. For a domain
  $R$ and a graph
  $\Gamma=(A;E)$ with the set of vertices
  $A$ and the set of edges
  $E$ denote by
  $\mathcal{L}_R(A;\Gamma)$,
  $\mathcal{M}_R (A;\Gamma)$, and
  $\mathcal{N}_{m,R}(A;\Gamma)$ the partially commutative
  $R$-algebra, partially commutative metabelian Lie
  $R$-algebra, and partially commutative nilpotent
  $R$-algebra of nilpotency degree
  $m$ respectively.

  Let
  $[u]$ be a Lie monomial in a (finite or infinite) set of generators
  $A=\{a_1,a_2,\dots\}$. \emph{Multi-degree} of
  $[u]$ is the vector
  $\overline{\delta}=(\delta_{1},\delta_2,\dots,)$, where
  $\delta_i$ is the number of occurrences of
  $a_i$ in
  $[u]$.

  For a Lie monomial
  $[u]$ of multi-degree
  $(\delta_{1},\delta_2,\dots,)$ put
  $\supp([u])=\{a_i\,|\,\delta_i\neq 0\}$. Extend this notation to the set of all Lie
  polynomials as follows. If
  $g=\sum_{j}\alpha_j [u_j]$ is a Lie polynomial then
  $\supp(g)=\bigcup_j \supp([u_j])$.

\subsection{Algebraic properties of partially commutative Lie algebras} \label{ssec21}
\subsubsection*{Isomorphisms}
  As far as we know partially commutative algebras are explored since
  80th. We start with the result for partially commutative associative algebras obtained
  by K.\,H.\,Kim, L.\,Makar-Limanov, J.\,Neggers, and F.\,W.\,Roush
\cite{KMNR80}. This result has been already mentioned in Subsection~%
 \ref{ssec12}. Nevertheless, we decided to give it in more details because this is one of
  the significant results for partially commutative algebras.

  For an undirected graph without loops
  $\Gamma= (A;E)$ and a domain
  $R$ denote by
  $\mathcal{A}_{R}(A;\Gamma)$ the partially commutative associative
  $R$-algebra defined by
  $\Gamma$.
  \begin{ttt}
    Let
    $\Gamma=(A;E)$ and
    $\Delta=(B;F)$ be undirected graphs without loops and
    $\mathbb{F}$ a field. The partially commutative
    associative algebras
    $\mathcal{A}_{\mathbb{F}}(A;\Gamma)$ and
    $\mathcal{A}_{\mathbb{F}}(B;\Delta)$ are isomorphic if and only if the graphs
    $\Gamma$ and
    $\Delta$ are isomorphic.
  \end{ttt}
  G.\,Duchamp and D.\,Krob in
\cite{DK93} generalized this result to the case of associative
  $R$-algebras where
  $R$ is an arbitrary domain. Besides, in the same paper they stated a criterium of existing isomorphism for partially commutative Lie algebras. The analogous criterium holds
  also for partially commutative metabelian Lie algebras~%
\cite{Por17}. In all cases, two algebras are isomorphic if and only if their defining graphs
  are isomorphic.

\subsubsection*{Bases} Finding linear bases is a significant problem because
  a linear basis is a very important tool for studying algebras.

  The first result on bases of partially commutative Lie algebras was obtained by G.\,Duchamp
  and D.\,Krob in
\cite{DK92}, but they did not give an explicit description of a  basis. Their algorithm was
  recursive. More precisely, let
  $R$ be a unital commutative ring and
  $\Gamma (A;E)$ a graph without loops. The corresponding partially commutative Lie
  $R$-algebra
  $\mathcal{L}_R(A;\Gamma)$ is considered, a totally disconnected set
  $B \subseteq A$ is chosen, and the problem was reduced to finding a linear basis of the
  algebra
  $\mathcal{L}_R(A\backslash B; \Gamma(A \backslash B))$.

  An explicit construction for bases of partially commutative Lie algebras was obtained in
\cite{Por11}. To make this description let us first recall a definition of
  Lyndon--Shirshov words.

  Denote by
  $A^{\star}$ the set of all associative and non-associative words (associative monomials with no coefficient) in
  $A$ respectively. We define the empty word by 1.

  Let us extend an arbitrary well order on
  $A$ to the lexicographic order on
  $A^{\star}$.

  An associative word
  $u$ is called an \emph{associative Lyndon--Shirshov word} if for any pair of nonempty words
  $v$ and
  $w$ such that
  $u=vw$ we have
  $wv<u$.

  A non-associative word (non-associative monomial with no coefficient)
  $[u]$ is called a \emph{Lyndon--Shirshov word} if
   \begin{enumerate}
     \item \label{ass}
       The word
       $u$ obtained from
       $[u]$ by omitting brackets is an associative Lyndon--Shirshov
       word;
     \item
       if
       $[u]=[[u_1],[u_2]]$, then
       $[u_1]$ and
       $[u_2]$ are Lyndon--Shirshov words (it follows from (1) that
       $u_1>u_2$);
     \item
       if
       $[u]=\bigl[[[u_{11}],[u_{12}]],[u_2]\bigr]$, then
       $u_2\geqslant u_{12}$.
   \end{enumerate}
   Denote the sets of all non-associative Lyndon--Shirshov words in
   $A$ by
   $\mathrm{LS}(A)$. It was shown in
\cite{Sh58} that the set
   $\mathrm{LS}(A)$ is a basis of free Lie
   $R$-algebra.

   For a partially commutative Lie algebra
   $\mathcal{L}_R(A;\Gamma)$ over a domain
   $R$ define  by induction \emph{partially commutaitve Lyndon--Shirshow
   words} (\emph{PCLS-words} for short).
   \begin{enumerate}
     \item
       All elements of
       $A$ are PCLS-words.
     \item
       a Lyndon--Shirshov word
       $[u]$ such that
       $\ell([u])>1$ is a PCLS-word if
       $[u]=[[v],[w]]$, where
       $[v]$ and
       $[w]$ are PCLS-words and there is an element in
       $\supp([v])$ such that it is not connected in
       $\Gamma$ with the first letter of
       $[w]$.
     \item
       There are no other PCLS-words.
   \end{enumerate}

   Denote the set of all PCLS-words of a partially commutative Lie
   $R$-algebra
   $\mathcal{L}_R(A;\Gamma)$ by
   $PCLS(A;\Gamma)$. Using the method of Gr\"{o}bner--Shirshov bases the explicit
   description of bases of partially commutative Lie algebras was obtained

   \begin{ttt}\label{lbasis}
\cite{Por11} Let
     $R$ be a unital commutative ring,
     $\Gamma$ a finite undirected graph without loops, and
     $A$ the set of vertices of
     $\Gamma$. Then the set
     $PCLS(A;\Gamma)$ is a basis of the partially commutative Lie
     $R$-algebra
     $\mathcal{L}_{R}(A;\Gamma)$.
   \end{ttt}

   A linear basis for a partially commutative nilpotent algebra can be easily
   obtained from a linear basis of a partially commutative algebra.
   \begin{ttt}\label{nbasis}
\cite{Por11-1} Let
     $R$ be a unital commutative ring,
     $\Gamma$ a finite undirected graph without loops, and
     $A$ the set of vertices of
     $\Gamma$. Then a basis of the partially commutative nilpotent
     $R$-algebra
     $\mathcal{N}_{R,m}(A;\Gamma)$ consists of all elements of
     $PCLS(A;\Gamma)$ whose lengthes are not greater than
     $m$.
   \end{ttt}

   In
\cite{Por11} and
 \cite{Por11-1} the set
   $A$ is supposed to be finite, but it is easy to see that this restriction is not
   essential and so, Theorem~%
\ref{lbasis} and Theorem~%
 \ref{nbasis} hold for algebras defined by infinite graphs
   $\Gamma$ as well.

   The problem of finding a linear basis for a partially commutative metabelian Lie
   algebra is also rather interesting. The explicit description of such basis was obtained
   in
\cite{PT13}. The idea used for constructing a basis for partially commutative metabelian
   Lie algebra is similar to one for a partially commutative Lie algebra. Namely, a basis
   is constructed by choosing some elements from a linear basis of a free metabelian
   Lie algebra of the corresponding variety.

   The basis for a free metabelian Lie algebra was obtained independently by
   L.\,A.\,Bokut
\cite{Bo63} and A.\,L.\,Shmelkin
 \cite{Sh64}. Let
   $A=\{a_1,a_2,\dots, a_n\}$ be the set of generators of free metabelian algebra. Then
   a basis of this algebra consists of the elements of the form
   $[a_{i_1},a_{i_2}, \dots, a_{i_k}]$, where
   $a_{i_1}>a_{i_2}$,
   $a_{i_2} \leqslant a_{i_3} \leqslant \dots \leqslant a_{i_k}$. Denote this set by
   $\mathrm{Bas}(A)$.

   Fix an arbitrary multi-degree
   $\overline{\delta}=(\delta_1,\delta_2, \dots, \delta_n)$, where
   $n=|A|$. Let
   $N=\sum_{i=1}^n \delta_i$,
   $A_{\overline{\delta}}=\{a_i\in A\,|\,\delta_i\neq 0\}$, and
   $b$ the smallest element of
   $A_{\overline{\delta}}$. Denote the connected components of the graph
   $\Gamma(A_{\overline{\delta}})$ by
   $\Delta_0,\Delta_1,\dots, \Delta_k$ in such a way that
   $b \in A(\Delta_0)$. Let
   $[u_{i}]\in\mathrm{Bas}(A)$ be an element of multi-degree
   $\overline{\delta}$ such that
   $[u_i]=[a_{j_{i,1}},b, a_{j_{i,3}}, \dots, a_{j_{i,N}}]$, where
   $a_{j_{i,1}}$ is the largest element of
   $A(\Delta_i)$. Denote by
   $B_{\overline{\delta}}(A;\Gamma)$ the subset
   $\{[u_1],[u_2], \dots, [u_k]\}$. Finally, put
   $$\mathrm{Bas}(A;\Gamma)=\bigcup_{\overline{\delta}}B_{\overline{\delta}}(A;\Gamma),$$
   where the union is taken on all multi-degrees.
   \begin{ttt} \label{mbasis}
{\upshape \cite{PT13}}
     Let
     $R$ be a unital commutative ring and
     $\Gamma$ a finite undirected graph without loops. Then the set
     $\mathrm{Bas}(A;\Gamma)$ is a basis of partially commutative metabelian
     Lie
     $R$-algebra
     $\mathcal{M}_{R}(A;\Gamma)$.
   \end{ttt}

   Note that Theorem~%
\ref{mbasis} also holds for infinitely generated partially commutative metabelian Lie algebras.

  \subsubsection*{Annihilators}
   Let
   $R$ be an infinite integral domain and
   $R[A]$ be the set of all commutative associative polynomials over
   $R$. The derived subalgebra
   $\mathcal{M}_R'(A;\Gamma)$ of the
   $R$-algebra
   $\mathcal{M}_R(A;\Gamma)$ is an
   $R[A]$-module with respect to the adjoint representation.

   Define the ideal
   $I^{\Gamma}_{i,j}$ of
   $R[A]$ as follows. If
   $a_i$ and
   $a_j$ are vertices belonging to different connected components in
   $\Gamma$ then put
   $I^{\Gamma}_{i,j}=0$. Suppose this is not so. Then for each path
   $(a_i,b_1,b_2,\dots, b_s,a_j)$ connecting these vertices in
   $\Gamma$ consider the associative monomial
   $b_1b_2\dots b_s$. Define
   $I^{\Gamma}_{i,j}$ as the ideal generated by all such monomials.
   \begin{ttt}\label{ann}
     {\upshape%
\cite{PT13}} Let
     $R$ be an infinite domain and
     $\Gamma=(A;E)$ a finite undirected graph without loops. For
     $a_i,a_j\in A$ if
     $a_i$ and
     $a_j$ are not adjacent in
     $\Gamma$ then the annihilator of
     $[a_i,a_j]$ in
     $\mathcal{M}_R(A;\Gamma)$ is equal to
     $I^{\Gamma}_{i,j}$.
   \end{ttt}

 \subsubsection*{Centralizers}
  As well as for isomorphisms, the first results for centralizers of partially commutative 
  algebras were obtained for associative ones. In 1980, K.\,H.\,Kim and F.\,W.\,Roush
  obtained a description of centralizers of monomials
\cite{KR80}.
  \begin{ttt}
     Let
     $R$ be a unital commutative ring,
     $\Gamma$ a finite undirected graph without loops, and
     $A$ the set of vertices of
     $\Gamma$. Let also
     $u$ be a monomial of degree
     $>0$ in the partially commutative associative
     $R$-algebra
     $\mathcal{A}_R(A;\Gamma)$ and let
     $v$ be a monomial of degree
     $>0$ in
     $\mathcal{A}_R(A;\Gamma)$ such that this monomial commutes with
     $u$. Finally, let
     $\Delta_1, \Delta_2, \dots ,\Delta_P$ be the connected components of
     $\Gamma^c(\supp(u))$. Write
     $u = u_1u_2\dots u_p$, where
     $\supp(u_i)=A(\Delta_i)$. Then
     $v$ is a product of generators
     $a$, such that
     $a\not \in \supp(u)$ but
     $a\leftrightarrow \supp(u)$, and words
     $w$ such that some power of
     $w$ equals one of the
     $u_i$.
  \end{ttt}

  It seems that the requirement of finiteness of the defining graph can be eliminated.

  We use the following notation
  Let
  $R$ be a domain and
  $L$ a Lie
  $R$-algebra. For
  $f,g\in L\backslash \{0\}$ we write
  $f \backsim g$ if
  $\alpha f=\beta g$ for some
  $\alpha, \beta \in R$. For any
  $g \in L$ the centralizer of
  $g$ is denoted by
  $C(g)$. We also put
  $\mathcal{C}(g)=C(g) \cap L'$.

  Unlike the case of partially commutative associative algebras, centralizers of elements
  of partially commutative Lie algebras over domains were described completely, i.e. for an
  arbitrary domain
  $R$ an explicit description for centralizers of all elements in any
  $R$-algebra
  $\mathcal{L}_R(A;\Gamma)$ was obtained.

  \begin{ttt}\label{gencase}
{\upshape \cite{Por12}}
    Let
    $R$ be a domain
    $\Gamma$ a finite undirected graph without loops, and
    $A$ the set of vertices of
    $\Gamma$. For an arbitrary element
    $g$ of the  Lie
    $R$-algebra
    $\mathcal{L}_R(A;\Gamma)$ denote by
    $\Delta_1,\Delta_2, \dots, \Delta_p$ all connected components of the graph
    $\Gamma^{c}(\supp(g))$. Then
    $g=\sum_{i=1}^p g_i$, where
    $\supp(g_i)= A(\Delta_i)$ for all
    $i=1,2,\dots p$, and
    $C(g)$ consists of elements of the form
    $h=\sum_{i=1}^p h_i+h^{(1)}$, where for each
    $i=1,2,\dots p$ either
    $h_i=0$ or
    $g_i\backsim h_i$. Moreover,
    $\supp(g)\leftrightarrow \supp(h^{(1)})$.
  \end{ttt}

  For partially commutative metabelian Lie algebras there is no complete description of centralizers. Nevertheless, in
\cite{PT13, Por15} some specific results were obtained.

  For
  $f\in \mathcal{M}_R'(A;\Gamma)$ and
  $g\in R[A]$ denote by
  $f.g$ the image of
  $f$ via the adjoint action by
  $g$.

  \begin{ttt}\label{centrmetgen}
    Let
    $R$ be a domain,
    $\Gamma$ a finite undirected graph without loops, and
    $A=\{a_1,a_2,\dots,a_n\}$ the set of vertices of this graph. Then for
    partially commutative metabelian Lie
    $R$-algebra
    $\mathcal{M}_R(A;\Gamma)$ the following statements hold.\\
    1) If
       $a_n$ is an isolated vertex in
       $\Gamma$ then
       $C(a_n)$ consists of the elements
       $v$ of the form
       \begin{equation*}
         v=\alpha_n a_n,
       \end{equation*}
       where
       $\alpha_a \in R$.\\
    2) If the degree of
       $x_n$ is equal to
       $1$ in
       $G$ (say, it is adjacent to
       $a_{n-1}$) then
       $C(a_n)$ consists of all elements
       $v$ of the form
       \begin{equation*}
         v=\alpha_{n-1}a_{n-1}+\alpha_n a_n,
       \end{equation*}
       where
       $\alpha_{n-1},\alpha_n\in R$.\\
    3) If
       $a_n$ is adjacent to
       $a_{r+1},\dots, a_{n-1}$ in
       $\Gamma$ ($r \leqslant n-3$), then
       $C(a_n)$ consists of all elements
       $v$ of the form
       \begin{equation*}
         v=\sum_{k=r+1}^n \alpha_k a_k +\sum_{r+1\leqslant i<j\leqslant n-1}[a_i,a_j].f_{ij},
       \end{equation*}
       where
       $\alpha_k \in R$,
       $f_{ij}\in R[A\backslash \{a_n\}]$.
  \end{ttt}

  There are some results on ``centralizers in the commutant''
  $\mathcal{C}(g)$ in partially commutative metabelian Lie
  $R$-algebras.

  \begin{ttt}\label{centinter}
    {\upshape%
\cite{PT13}}
    Let
    $R$ be a domain,
    $\Gamma$ a finite undirected graph without loops, and
    $A=\{a_1,a_2,\dots,a_n\}$ the set of vertices of this graph.
    Then in the partially commutative metabelian Lie
    $R$-algebra
    $\mathcal{M}_R(A;\Gamma)$ the following equation holds.
    $$\mathcal{C}\bigl(\sum_{j=1}^m \alpha_{i_j} a_{i_j}\bigr)=\bigcap_{j=1}^{m}
    \mathcal{C}(a_{i_j})$$
    for any elements
    $a_{i_1},a_{i_2},\dots, a_{i_m}$ and for any
    $\alpha_{i_1},\alpha_{i_2},\dots,\alpha_{i_m}\in R\backslash\{0\}$.
  \end{ttt}

  Theorem~%
\ref{centinter} has some corollaries for partially commutative metabelian Lie
  $R$-algebras defined by specific graphs.

  \begin{ccc} \label{centralizators}%
    {\upshape%
\cite{Por15}}
    Let
    $R$ be a domain,
    $C_n$ a cycle of length
    $n\geqslant 3$, and
    $A=\{a_1,a_2,\dots, a_n\}$ the set of vertices of
    $C_n$. Then the following statements hold in
    $\mathcal{M}_R(A;C_n)$.\\
    a) if
    $a_i$ and
    $a_j$ are adjacent then
    $\mathcal{C}(\alpha a_i+\beta a_n)=0$ for any
    $\alpha,\beta \in R\backslash\{0\}$.\\
    b) If
    $a_i$ and
    $a_j$ are not adjacent then for any
    $\alpha,\beta \in R\backslash\{0\}$ the set
    $\mathcal{C}(\alpha a_i+\beta a_j)$ consists of linear combinations of non-zero Lie
    monomials
    $[u_r]$ such that
    $A(\supp([u_r]))=A\backslash\{a_i,a_j\}$. Moreover, any element of
    $\mathcal{C}(\alpha x_i+\beta x_j)$ can be represented in the form
    $f=[a_{i-1},a_{i+1}].g$ for some
    $g \in R[A]$. \\
    c) If
    $m\geqslant 3$, then
    $\mathcal{C}(\sum_{j=1}^m \alpha_j x_{i_j})=0$ for any
    $a_{i_1},a_{i_2}, \dots, a_{i_m}$ and
    $\alpha_1,\alpha_2,\dots, \alpha_m \in R\backslash \{0\}$.
   \end{ccc}

   \begin{ccc}\label{centlincomb}%
    {\upshape%
\cite{PT13}} Let
    $R$ be a domain,
    $\Gamma$ a finite tree, and
    $A=\{a_1,a_2,\dots,a_n\}$ the set of vertices of this tree. Suppose that
    $g= \sum_{j=1}^{m}\alpha_j a_{i_j}$, where
    $m\geqslant 2$ and
    $\alpha_j\in R\backslash\{0\}$ for
    $j=1,2, \dots, m$. Then
    $\mathcal{C}(g)=0$ in
    $\mathcal{M}_R(A;\Gamma)$.
  \end{ccc}

 \subsection{Logical properties of partially commutative and partially commutative
    metabelian Lie algebras} \label{ssec22}
 \subsubsection*{Universal equivalence}

   Conditions of universal equivalence of partially commutative Lie
   $R$-algebras over domains
   $R$ were studied in the series of papers
\cite{PT13,Por15,Por17-1,Por17-2,Por20}. The problem of finding a criteria of universal
   equivalence on the entire class of partially commutative (metabelian) Lie
   $R$-algebras
   seems to be very complicated. So, this problem is considered on some specific classes
   of
   $R$-algebras. Although the methods used in partially commutative and partially
   commutative metabelian Lie
   $R$-algebras differ essentially the results turned out to be
   similar. The first criteria of universal equivalence were obtained for Lie
   $R$-algebras defined by cycles and trees.

   \begin{ttt}\label{pccycunieq}
{\upshape \cite{Por17-1}}
     Let
     $R$ be a domain and
     $C_n=(A;E)$ and
     $C_m= (B;F)$ cyclic graphs such that
     $|A|=n$,
     $|B|=m$ with
     $n,m \geqslant 3$. Then the partially commutative Lie
     $R$-algebras
     $\mathcal{L}_R(A;C_n)$ and
     $\mathcal{L}_R(B;C_m)$ are universally equivalent if and only if
     $n=m$.
   \end{ttt}

   The analogous result holds also for partially commutative metabelian Lie
   $R$-algebras with some restriction on a domain
   $R$.

   \begin{ttt}\label{pcmetcycunieq}
{\upshape \cite{Por15}}
     Let
     $R$ be a domain containing
     $\mathbb{Z}$ as a subring and
     $C_n=(A;E)$ and
     $C_m= (B;F)$ cyclic graphs such that
     $|A|=n$,
     $|B|=m$ with
     $n,m \geqslant 3$. Then the partially commutative metabelian Lie
     $R$-algebras
     $\mathcal{M}_R(A;C_n)$ and
     $\mathcal{M}_R(B;C_m)$ are universally equivalent if and only if
     $n=m$.
   \end{ttt}

   A criterium of universal equivalence of partially commutative
   $R$-algebras, where
   $R$ is a domain was found in
\cite{Por17-1}. The analogous result in the metabelian case was obtained in
 \cite{PT13}. Despite that result was obtained only for partially commutative
  metabelian rings it can be easily generalized to the case of partially commutative
  metabelian
  $R$-algebra, where
  $R$ is a domain containing
  $\mathbb{Z}$ as a subring.

  It turned out that these results can be generalized to the case
   of algebras defined by graphs with at most countably many vertices.

  Let
  $\Delta=(A;E)$ be a graph. Denote by
  $A^*$ the set obtained from
  $A$ by deleting all end-points of
  $\Delta$ and put
  $\Delta^*=\Delta(A^*)$. Actually, the definition of the graph
  $\Delta^*$ coincides with one given in Theorem~%
\ref{abelguni}.

  We say that a tree
  $\Gamma$ (finite or infinite) is a \emph{tree of finite type} if the tree
  $\Gamma^*$ is finite and  a \emph{tree of infinite type} if
  $\Gamma^*$ is infinite.

  Graphs
  $\Gamma$ and
  $\Delta$ are \emph{mutually locally embeddable} if any finite subgraph of each graph
  $\Gamma$ and
  $\Delta$ is isomorphically embeddable to the other one. The following theorem shows
  that two partially commutative Lie algebras are generated by graphs of different types
  then this algebras are not universally equivalent.

   \begin{ttt}\label{pcafin-infin}
{\upshape \cite{Por17-2}} Let
     $R$ be a domain,
     $\Gamma=(A;E)$ a tree of infinite type, and
     $\Delta=(B;F)$ a tree of finite type, where the sets
     $A$ and
     $B$ are at most countable. Then the partially commutative Lie
     $R$-algebras
     $\mathcal{L}_R(A;\Gamma)$ and
     $\mathcal{L}_R(B;\Delta)$ are not universally equivalent for
     any domain
     $R$.
  \end{ttt}
  The following theorem is an analogue of Theorem~%
\ref{pcafin-infin} for partially commutative metabelian Lie algebras.
  \begin{ttt}
{\upshape \cite{Por17-2}} Let
     $R$ be a domain containing
     $\mathbb{Z}$ as a subring,
     $\Gamma=(A;E)$ a tree of infinite type, and
     $\Delta=(B;F)$ a tree of finite type, where the sets
     $A$ and
     $B$ are at most countable. Then the partially commutative metabelian Lie
     $R$-algebras
     $\mathcal{M}_R(A;\Gamma)$ and
     $\mathcal{M}_R(B;\Delta)$ are not universally equivalent.
   \end{ttt}

   So, there are separate criteria for Lie algebras defined by graphs of finite type
   and for Lie algebras defined by graphs of infinite type. The case of finite graphs is
   considered in Theorem~%
\ref{univeqfintypelie} for partially commutative Lie algebras and in Theorem~%
  \ref{univeqfintypemetlie} for partially commutative metabelian Lie algebras.
   \begin{ttt}\label{univeqfintypelie}
{\upshape \cite{Por17-1,Por17-2}} Let
     $R$ be a domain,
     $\Gamma=(A;E)$ and
     $\Delta=(B;F)$ be trees of finite type such that
     $A$ and
     $B$ are at most countable and
     $|A^*|\geqslant 2$,
     $|B^*|\geqslant 2$. Then partially commutative Lie
     $R$-algebras
     $\mathcal{L}_R(A;\Gamma)$ and
     $\mathcal{L}_R(B;\Delta)$ are universally equivalent if and only if
     $\Gamma^*\simeq \Delta^*$.
  \end{ttt}

  \begin{ttt} \label{univeqfintypemetlie}
     {\upshape \cite{PT13,Por17-2}} Let
     $R$ be a domain containing
     $\mathbb{Z}$ as a subring,
     $\Gamma=(A;E)$ and
     $\Delta=(B;F)$ be trees of finite type such that
     $A$ and
     $B$ are at most countable and
     $|A^*|\geqslant 2$,
     $|B^*|\geqslant 2$. Then the partially commutative metabelian Lie
     $R$-algebras
     $\mathcal{L}_R(A;\Gamma)$ and
     $\mathcal{L}_R(B;\Delta)$ are universally equivalent if and only if
     $\Gamma^*\simeq \Delta^*$.
   \end{ttt}

   Finally, the following two theorems provide criteria of universal equivalence of
   partially commutative and partially commutative metabelian Lie algebras generated by
   graphs of infinite type.
   \begin{ttt}\label{univeqinfintypelie}
{\upshape \cite{Por17-2}}  Let
    $R$ be a domain,
    $\Gamma=(A;E)$ and
    $\Delta=(B;F)$ trees of infinite type and these trees have at most countable sets
    of vertices. Then the partially commutative Lie
    $R$-algebras
    $\mathcal{L}_R(A;\Gamma)$ and
    $\mathcal{L}_R(B;\Delta)$ are universally equivalent if and only if
    $\Gamma^*$ and
    $\Delta^*$ are mutually locally embeddable.
  \end{ttt}
  \begin{ttt} \label{univeqinfintypeliemetab}
{\upshape \cite{Por17-2}}  Let
    $R$ be a domain containing
    $\mathbb{Z}$ as a subring,
    $\Gamma=(A;E)$ and
    $\Delta=(B;F)$ trees of infinite type and these trees have at most countable sets
    of vertices. Then the partially commutative metabelian Lie
    $R$-algebras
    $\mathcal{M}_R(A;\Gamma)$ and
    $\mathcal{M}_R(B;\Delta)$ are universally equivalent if and only if
    $\Gamma^*$ and
    $\Delta^*$ are mutually locally embeddable.
   \end{ttt}

   Note that the statements similar to ones in Theorems
\ref{pcafin-infin}--%
 \ref{univeqinfintypeliemetab} were obtained for countably generated partially commutative
   metabelian groups. The following theorem shows that no partially commutative metabelian
   group defined by a graph of finite type can be universally equivalent to one defined by
   a graph of infinite type.

   \begin{ttt}\label{grfin-infin}
{\upshape \cite{Por17-2}}   Let
     $\Gamma=(A;E)$ a tree of infinite type, and
     $\Delta=(B;F)$ a tree of finite type, where the sets
     $A$ and
     $B$ are at most countable.Then the partially commutative metabelian groups
     $G=(\mathfrak{A}^2,\Gamma)$ and
     $H=(\mathfrak{A}^2, \Delta)$ are not universally equivalent.
   \end{ttt}

   The following two theorems establish criteria of universal equivalence of partially
   commutative metabelian groups defined by graphs of finite (Theorem~%
\ref{univeqfintypegr}) and infinite (Theorem~%
 \ref{univeqinfintypegr}) type.

   \begin{ttt}\label{univeqfintypegr}
{\upshape \cite{Por17-2}} Let
     $\Gamma=(A;E)$ and
     $\Delta=(B;F)$  trees of finite type such that
     $A$ and
     $B$ are at most countable, at least one of them is countable, and
     $|A^*|\geqslant 2$,
     $|B^*|\geqslant 2$. Then the partially commutative metabelian groups
     $G=(\mathfrak{A}^2,\Gamma)$ and
     $H=(\mathfrak{A}^2, \Delta)$ are universally equivalent if and only if
     $\Gamma^*\simeq \Delta^*$.
   \end{ttt}

   Actually, Theorem~
\ref{univeqfintypegr} generalizes the criterium of universal equivalence of partially commutative
   metabelian groups defined by finite trees
 \cite{GT11-1}.

   \begin{ttt}\label{univeqinfintypegr}
{\upshape \cite{Por17-2}}  Let
    $R$ be a domain,
    $\Gamma=(A;E)$ and
    $\Delta=(B;F)$ trees of infinite type and these trees have  at most countable sets of vertices. Then the partially commutative metabelian groups
    $G=(\mathfrak{A}^2,\Gamma)$ and
    $H=(\mathfrak{A}^2, \Delta)$ are universally equivalent if and only if
    $\Gamma^*$ and
    $\Delta^*$ are mutually locally embeddable.
  \end{ttt}

   It is rather easy to see that the condition on cardinalities of sets
   $A$ and
   $B$ in Theorems~%
\ref{pcafin-infin}--%
 \ref{univeqinfintypegr} can be excluded.

   The results in
\cite{Por17-1,Por17-2} for partially commutative Lie algebras were generalized in
 \cite{Por20} as follows.

   \begin{ttt}
    Let
    $R$ be a domain,
    $\Gamma=(A;E)$ and
    $\Delta=(B;F)$ finite undirected graphs without loops, triangles, squares, and isolated
    vertices. Then the partially commutative Lie
    $R$-algebras
    $\mathcal{L}_R(A;\Gamma)$ and
    $\mathcal{L}_R(B;\Delta)$ are universally equivalent if and only if
    $\Gamma^*\simeq \Delta^*$  and quantities of two-vertex connected components in
    $\Gamma$ and
    $\Delta$ are equal.
   \end{ttt}

   In
\cite{Por17-1}, it was shown that the class of partially commutative Lie algebras
   defined by finite trees is not
   distinguished in the class of all finitely generated partially commutative algebras by
   universal theories.

 \subsubsection*{Elementary equivalence}
   The problem of finding criteria of elementary equivalence for partially commutative
   and partially commutative metabelian Lie algebras was studied in
\cite{Por17}. This problem can be considered for algebras not only in ``classical''
   signature, but also when Lie algebras are considered as two-sorted algebraic systems.
   Namely, there are three operations considered: addition and multiplication of elements
   in the algebra and multiplication of an element in the basic field by an element in
   the algebra.

   Criteria of elementary equivalence of partially commutative and partially commutative
   metabelian Lie algebras were found in the case when Lie algebras over a field are
   considered as two-sorted systems and for Lie rings.

   \begin{ttt}\label{elemeqtwobase}
{\upshape \cite{Por17}}
     \renewcommand{\theenumi}{{\upshape \arabic{enumi}}}
       Let
       $\mathbb{F}$ be a field and
       $\Gamma=(A;E)$ and
       $\Delta=(B;F)$ finite undirected graphs without loops.
       \begin{enumerate}
         \item
           The partially commutative Lie algebras
           $\mathcal{L}_{\mathbb{F}}(A;\Gamma)$ and
           $\mathcal{L}_{\mathbb{F}}(B;\Delta)$, considered as two-sorted algebraic systems
           are elementary equivalent if and only if
           $\Gamma \simeq \Delta$.
         \item
           The partially commutative metabelian Lie algebras
           $\mathcal{M}_{\mathbb{F}}(A;\Gamma)$ and
           $\mathcal{M}_{\mathbb{F}}(B;\Delta)$, considered as two-sorted algebraic systems
           are elementary equivalent if and only if
           $\Gamma \simeq \Delta$.
       \end{enumerate}
   \end{ttt}

   \begin{ttt}\label{elemeqfinring}
 {\upshape \cite{Por17}}
     Let
     $\Gamma=(A;E)$ and
     $\Delta=(B;F)$ finite undirected graphs without loops.
     \renewcommand{\theenumi}{{\upshape \arabic{enumi}}}
     \begin{enumerate}
       \item
         The partially commutative Lie rings
         $\mathcal{L}_{\mathbb{Z}}(A;\Gamma)$ and
         $\mathcal{L}_{\mathbb{Z}}(B;\Delta)$ are elementary equivalent if and only if
         $\Gamma\simeq \Delta$.
       \item
         The partially commutative metabelian Lie rings
         $\mathcal{M}_{\mathbb{Z}}(A;\Gamma)$ and
         $\mathcal{M}_{\mathbb{Z}}(B;\Delta)$ are elementary equivalent if and only if
         $\Gamma\simeq \Delta$.
     \end{enumerate}
 \end{ttt}

 \end{document}